\documentclass[a4paper,12pt]{article}
\usepackage[T1]{fontenc}
\usepackage[utf8]{inputenc}
\usepackage{lmodern}
\usepackage{amsmath}
\usepackage{amsfonts}
\usepackage{amssymb}
\usepackage[usenames,dvipsnames,svgnames,table]{xcolor}
\usepackage{graphicx}
\usepackage{amsthm}
\usepackage{tikz}
\usepackage{microtype}
\usepackage{siunitx}
\usetikzlibrary{angles, quotes}
\usepackage{mathrsfs}
\usepackage[hmargin=2cm,vmargin=3cm]{geometry}
\usepackage{enumitem}
\usepackage{bm}
\usepackage{algorithm}
\usepackage{algorithmic}
\usepackage{booktabs,multirow}
\usepackage{colortbl}
\usepackage{diagbox}
\usepackage{comment}
\usepackage{bbold}
\usepackage{pgfplots}
\usepackage{subcaption}
\usepackage{listings}
\usepackage{appendix}
\usepackage{calc}

\theoremstyle{definition}

\theoremstyle{plain}

\theoremstyle{remark}
\newtheorem{remark}{Remark}[section]

\providecommand{\keywords}[1]
{
	\small	
	\textbf{\textit{Keywords ---}} #1
}

\newcommand{\nug}{\bm{n}}
\newcommand{\mup}[1]{\mu_{{#1}} }
\newcommand{\discrd}[2]{\mathcal{#1}_{\texttt{#2}}}

\newcommand{\fespacebm}[1]{\bm{\mathcal{#1}}}
\newcommand{\enadm}[1]{\bm{\mathcal{#1}}}

\newcommand{\rfespacebm}[2]{\bm{#1}_{{#2 }}}

\newcommand{\norm}[2]{\|#1\|_{{#2}}}

\newcommand{\set}[3]{\{#1_{{#2}}\}_{{#2\in\{#3:\uppercase{#2}\}}}}
\newcommand{\setbbb}[3]{\{#1_{{#2}}\}_{{#2\in\{#3:\uppercase{#2}\}}}}

\newcommand{\setbis}[3]{\{#1\}_{{#2\in\{#3:\uppercase{#2}\}}}}

\newcommand{\matr}[2]{\big[#1\cdots #2\big]}
\newcommand{\redsum}[3]{\sum\limits_{#1 \in\{#2:\uppercase{#1}^{#3}\}}}

\newcommand{\Iparr}[1]{\big(#1\big)}

\newcommand{\taug}{\bm{\tau}}

\newcommand{\ie}{{\emph{i.e.}}}

\newcommand{\chf}{c_{{\rm{HF}}}}
\newcommand{\betahf}{{\beta}_{{\rm{HF}}}}

\newcommand{\proj}[2]{\Pi_{{#1}}^{#2}}

\newcommand{\netoi}{\mathbb{N}^{*}}

\newcommand{\nuf}{\nu_{\mathcal{F}}}
\newcommand{\omegaa}{\Omega}
\newcommand{\omegaahat}{\widehat{\Omega}}
\newcommand{\gammaa}{\Gamma}
\newcommand{\gammaahat}{\widehat{\Gamma}}
\newcommand{\gammaad}{\Gamma^{\textup{D}}}
\newcommand{\gammaac}{\Gamma^{\textup{c}}}
\newcommand{\gammaan}{\Gamma^{\textup{N}}}
\newcommand{\gammaachat}{\widehat{\Gamma}^{\textup{c}}}
\newcommand{\gammaadhat}{\widehat{\Gamma}^{\textup{D}}}
\newcommand{\gammaanhat}{\widehat{\Gamma}^{\textup{N}}}

\newcommand{\projconv}[2]{\big[#1\big]_{{#2}}}

\newcommand{\getfem}{\texttt{getfem}}

\newcommand{\prpu}{\delta_{\bm{u}}}

\newcommand{\rdp}{\mathbb{R}_{+}}

\newcommand{\kc}{k^{\rm{cv}}}

\newcommand{\varphibm}{\bm{\varphi}}

\newcommand{\e}{e}

\newcommand{\ubmka}[1]{{\bm{u}_{#1}}}

\newcommand{\Ubm}{\bm{U}}

\newcommand{\dimP}{\mathcal{N}^{\rm{HF}}}
\newcommand{\dimPred}{N}

\newcommand{\F}{F}

\newcommand{\Ubmka}[1]{{\bm{U}}_{#1}}

\newcommand{\xibm}{\bm{\xi}}
\newcommand{\Xibm}{\bm{\Xi}}

\newcommand{\setbbbbis}[3]{\{#1_{{#2}}\}_{{#2\in\{#3:\uppercase{#2}\}}}}

\newcommand{\ubmred}{{\bm{u}_{\dimPred}}}

\newcommand{\ubmredvka}[2]{\bm{u}_{\dimPred,#2}^{#1}}

\newcommand{\Vect}[5]{(#1(#2))_{{#3\in\{#4:#5\}}}}

\newcommand{\Ubmredka}[1]{\bm{U}_{\dimPred,#1}}

\newcommand{\ubmredka}[1]{{\bm{u}_{\dimPred,#1}}}

\newcommand{\Rd}[1]{\mathbb{R}^{#1}}

\newcommand{\Pgamn}{P^{\nug}_{\gamma,g}}
\newcommand{\Pgamnz}{P^{\nug}_{\gamma,0}}
\newcommand{\Pgamt}{P^{\taug}_{\gamma}}
\newcommand{\Agam}{a_{\gamma}^{\nug \taug} }
\newcommand{\Agamn}{a^{\nug}_{\gamma}}
\newcommand{\Resn}{r^{\nug}_{\gamma}}

\newcommand{\Res}{r^{\nug \taug}_{\gamma}}

\newcommand{\Bgamn}{b^{\nug}_{\gamma}}
\newcommand{\Bgamt}{b^{\taug}_{\gamma}}
\newcommand{\Bgam}{b^{\nug \taug}_{\gamma}}

\newcommand{\bspan}{{\bm{Span}}}

\newcommand{\vbm}{{\bm{v}}}
\newcommand{\ubm}{{\bm{u}}}
\newcommand{\deltaubm}{{\bm{\delta u}}}
\newcommand{\Deltaubm}{{\bm{\Delta U}}}
\newcommand{\deltaubmv}[1]{\bm{\delta u}^{#1}}

\newcommand{\Deltaubmka}[1]{{\bm{\Delta U}_{#1}}}
\newcommand{\wbm}{{\bm{w}}}
\newcommand{\xbm}{{\bm{x}}}
\newcommand{\ellbm}{{\bm{\ell}}}

\newcommand{\sigmabm}{{\bm{\sigma}}}
\newcommand{\vn}{{v_{\nug}}}
\newcommand{\un}{{u_{\nug}}}

\newcommand{\vt}{{\vbm_{\taug}}}
\newcommand{\ut}{{\ubm_{\taug}}}
\newcommand{\sigman}{{\sigma_{\nug\nug}}}
\newcommand{\sigmat}{{\sigmabm_{\nug\taug}}}
\newcommand{\zerobm}{{\bm{0}}}

\newcommand{\JN}{{J^{\textup{\rm{\tiny{Nitsche}}}}}}
\newcommand{\JNF}{{J_{\mathcal{F}}^{\textup{\rm{\tiny{Nitsche}}}}}}
\newcommand{\oph}{h}
\newcommand{\ophc}{\oph_{c}}
\newcommand{\ophd}{\oph_{d}}
\newcommand{\ophn}{\oph_{n}}

\newcommand{\Sum}[3]{  \sum\limits_{#1 \in\{#2:#3\} } }

\newcommand{\deltaubmredka}[1]{{\bm{\delta u}_{N,#1}}}

\newcommand{\Deltaubmredka}[1]{{\bm{\Delta U}_{N,#1}}}

\newcommand{\deltaubmredvka}[2]{{\bm{\delta u}_{N,#2}^{#1}}}

\newcommand{\Agamns}{A^{\nug}_{\gamma,s}}

\newcommand{\Bgamns}{B^{\nug}_{\gamma,s}}

\newcommand{\AgamnM}{A^{\nug}_{\gamma,N}}
\newcommand{\BgamnM}{B^{\nug}_{\gamma,N}}
\newcommand{\ResnM}{{R}^{\nug}_{\gamma,N}}

\newcommand{\Agamnm}{A^{\nug}_{\gamma}}
\newcommand{\Agamm}{A^{\nug\taug}_{\gamma}}
\newcommand{\Bgamnm}{B^{\nug}_{\gamma}}
\newcommand{\Bgamm}{B^{\nug\taug}_{\gamma}}
\newcommand{\Resnm}{{R}^{\nug}_{\gamma}}

\newcommand{\Resm}{{R}^{\nug\taug}_{\gamma}}

\newcommand{\Alphas}[2]{\alpha_{#1}^{#2}}

\newcommand{\AgamnMs}{A^{\nug}_{\gamma,N,s}}

\newcommand{\BgamnMs}{B^{\nug}_{\gamma,N,s}}

\newcommand{\Z}{Z}

\newcommand{\ubmv}[1]{\ubm^{#1}}

\newcommand{\erbpr}{e_{\dimPred,\dimDred}^{\ubm}}
\newcommand{\erbdu}{e_{\dimPred,\dimDred}^{\lambda}}
\newcommand{\eimtr}[1]{\mathcal{E}^{#1}_{\texttt{train}}}

\newcommand{\eimcar}[1]{{S^{#1}}}
\newcommand{\eimop}[1]{E^{#1}}
\newcommand{\eimopN}[1]{E^{#1}_{\dimPred}}
\newcommand{\is}[1]{n_{s}^{#1}}
\newcommand{\js}[1]{m_{s}^{#1}}
\newcommand{\im}[1]{n_{i}^{#1}}
\newcommand{\jm}[1]{m_{i}^{#1}}

\newcommand{\eimmatrUp}[1]{Q^{#1}}

\newcommand{\eimsec}[1]{T^{#1}}
\newcommand{\dimDred}{R}
\newcommand{\lambdared}{\lambda_{\dimDred}}

\newcommand{\tole}[1]{\delta_{\texttt{#1}}}
\newcommand{\indset}[2]{\{#1{:}#2\}}
\newcommand{\stiff}{\mathbb{W}}

\newcommand{\projconvb}[2]{[#1]_{{#2}}}

\newcommand{\Gdif}[1]{\mathbb{G}_{#1}}
\newcommand{\mulame}{\mu^{\rm{Lamé}}}
\newcommand{\setcalmuhf}[4]{\{#1_{{#2}}(#4)\}_{{#2\in\{#3: \mathcal{\uppercase{#2}}^{\rm{HF}} \}}}}

\newcommand{\snc}{s^{\rm{c}}}

\newcommand{\AlphaNNs}[2]{\alpha_{\dimPred,#1}^{#2}}
\newcommand{\AlphaNN}[1]{\alpha_{\dimPred}^{#1}}
\newcommand{\phibm}{\bm{\phi}}

\newcommand{\Fs}{F_{s}}
\newcommand{\Rhons}{{\Theta}^{\nug}_{\gamma,s}}
\newcommand{\Rhotm}{{\Theta}^{\taug}_{\gamma}}
\newcommand{\Rhonm}{{\Theta}^{\nug}_{\gamma}}
\newcommand{\Rhon}{\theta^{\nug}_{\gamma}}
\newcommand{\RhonMs}{{\Theta}^{\nug}_{\gamma,N,s}}
\newcommand{\FM}{F_{N}}
\newcommand{\FMs}{\F_{N,s}}
\newcommand{\Rhot}{\theta^{\taug}_{\gamma}}
\newcommand{\RhonM}{{\Theta}^{\nug}_{\gamma,N}}
\newcommand{\Bgamtm}{B^{\taug}_{\gamma}}

\newcommand{\erb}[1]{e_{\dimPred}^{#1}}
\newcommand{\thetan}{\theta^{\nug}}
\newcommand{\thetat}{\theta^{\taug}}

\makeatother
\usepackage[english]{babel}

\title{\textbf{ \centering A reduced basis method for frictional contact problems formulated with Nitsche's method}  }

\author{\textsf{  \small \sc  Idrissa Niakh \footnote{ EDF R\&D, 7 Boulevard Gaspard Monge, 91120 Palaiseau, France} \footnote{ CERMICS, École des Ponts, 6-8 avenue Blaise Pascal, 77455 Marne-la-Vallée cedex 2, France} \footnote{INRIA Paris, 2 Rue Simone Iff, 75012 Paris, France}\ , Guillaume DROUET\footnotemark[1],   Virginie Ehrlacher \footnotemark[2] \footnotemark[3]\ , Alexandre Ern \footnotemark[2] \footnotemark[3]\ }}		
\begin{document}
\maketitle
\begin{abstract}
We develop an efficient reduced basis method for the frictional contact problem formulated using Nitsche's method. We focus on the regime of small deformations and on Tresca friction. The key idea ensuring the computational efficiency of the method is to treat the nonlinearity resulting from the contact and friction conditions by means of the Empirical Interpolation Method.
The proposed algorithm is applied to the Hertz contact problem between two half-disks with parameter-dependent radius. We also highlight the benefits of the present approach with respect to the mixed (primal-dual) formulation.
\end{abstract}


\section{Introduction}\label{section:4:intor}
The reduced basis method (RBM) is a model reduction technique~\cite{prud2002reliable,Buffa-2012,quarteroni2015reduced,hesthaven2016certified}.  The goal is to reduce the complexity of a parametrized model problem in computational studies where the parameters vary. The idea is to replace the \itshape high fidelity \normalfont (HF) discretization space, which is supposed to be of very large dimension, by a small-dimensional subspace (called \itshape reduced \normalfont space) which can be constructed by sampling the HF model. This allows one to organize the calculations in two phases. The first phase, called offline, is the construction phase of the reduced model. For this purpose, one considers a sample of the parameter space (assumed to be sufficiently representative) for which expensive calculations are performed by solving the HF problem for each parameter of the sample. The results of these calculations are then used to construct a small-dimensional subspace of the HF space, and the reduced model is built by replacing the HF space by the reduced subspace in a Galerkin-type approximation of the model. The second phase, called online, is a phase in which a large number of new values of the parameter are considered, for which accurate approximations of the HF solution are calculated by using the reduced model. The online phase is where substantial computational gains are achieved.

In this work, we are interested in the application of the RBM to the contact problem formulated with Nitsche's method. We focus on the regime of small deformations and on Tresca friction. The problem of mechanical contact~\cite{Johnson-1987,Wriggers-2006} with or without friction is  present in many structural  problems encountered in several industrial fields. The variational formulation of this problem leads to a variational inequality of the first or second kind depending on whether there is friction or not~\cite{Fichera-1964, Duvaut-1972}. There are different approaches to solve these variational inequalities. We can mention mixed (primal-dual) methods~\cite{Fortin-1983,Kikuchi-1988, Baillet-2006} where  Lagrange multipliers are introduced to enforce the contact and friction conditions. In this case, the problem to be solved is a saddle-point problem where one seeks a primal unknown (the displacement) and a dual unknown (the contact forces). One of the difficulties with these methods is that they require the contact operator to satisfy an inf-sup condition. In the literature, there is already some work on model reduction for the frictionless contact problem in the framework of a mixed formulation. For example, \cite{Haasdonk-2012} derives model reduction methods in the general framework of variational inequalities including the unilateral contact problem. In~\cite{Balajewicz-2016}, a projection-based method is proposed to reduce the contact problem under small deformations. In~\cite{Fauque-2018}, an application of the hyper-reduction technique is presented for the contact problem under small deformations. We also mention~\cite{Benaceur-2020} where a new dual basis construction is proposed for the RBM applied to the unilateral contact problem under large deformations. Finally, the recent work~\cite{Niakh-2022} proposes a stable and efficient model reduction method for the unilateral contact problem. Let us also mention~\cite{zeka2022preliminary} where the authors use the Progressive Generalized Method in order to build reduced-order models for problems with multiple contacts, and \cite{LeBRFD:23} where the authors use a hyper-reduction approach based on a reduced integration domain for the dual reduced basis.

In contrast to the mixed formulation approach, there are other methods to approximate the mechanical contact problem which are purely primal, \ie, they do not require the introduction of additional unknowns. These methods have the advantage of leading to unconstrained minimization problems (thus easier to solve) but do not guarantee that the contact and friction conditions are strictly satisfied. One example are penalty methods~\cite{Sofonea-2012}.
Here, we focus on another primal approach based on Nitsche's method~\cite{Nitsche-1971}. This method was originally introduced for the reformulation of Dirichlet boundary conditions and extended in~\cite{Chouly-2013} to the frictionless contact problem in the framework of the finite element method. The main characteristic of Nitsche's method is that it is consistent, in contrast to classical penalty methods. In the last few years, many contributions have been made to this approach. For example, in~\cite{Chouly-2014}, the method is extended to the case of contact with Tresca's friction; in~\cite{Chouly-2015}, symmetric and nonsymmetric variants are presented; in~\cite{Mlika-2017}, an extension to Coulomb's friction and large deformations is discussed; in~\cite{Chouly-2020}, a nonconforming high-order discretization is considered; in~\cite{Chouly-2022}, existence results for the contact problem with Coulomb friction are given in the context of static and dynamic finite element formulations. A state of the art on recent advances on Nitsche's method can be found in~\cite{Chouly-2017}. To the best of our knowledge, there is no previous work on model reduction for the contact problem formulated with Nitsche's method.

In this paper, we propose to fill this gap for the frictionless contact problem and the contact problem with Tresca friction. The main challenge is that the classical RBM leads to an inefficient reduced model owing to the nonlinearity of Nitsche's formulation (even with small deformations). To overcome this problem, we propose a combination of the RBM with the Empirical Interpolation Method (EIM)~\cite{Barrault-2004,Maday-2009}. The realization of this idea is by no means straightforward since one needs to consider at the same time the parameter value and the iteration counter in the nonlinear iterative solver. The second important point addressed in this work is the comparison of the present approach with the inf-sup stable mixed formulation in terms of accuracy and efficiency. The two key advantages offered by Nitsche's method are the handling of unconstrained minimization problems and a higher effectiveness of the RBM since it is well-known that the dual basis is particularly hard to compress, as was highlighted in particular in~\cite{kollepara2023limitations}.

The rest of this paper is organized as follows. In Section~\ref{section:4:modpro}, we briefly recall the unilateral contact problem with friction and its variational formulation under the assumption of small deformations. In Section~\ref{section:4:nitmet}, we derive the formulation of this problem using Nitsche's method in a form suitable to the RBM. In Section~\ref{section:4:rbm}, we present our main result, namely the procedure for building the reduced model with Nitsche's method using the RBM and the EIM. In Section~\ref{section:4:numres}, we provide numerical results showcasing the efficiency and the robustness of the proposed procedure and comparing it to the mixed formulation. We consider as test case the Hertz contact problem between two half-disks with parameter-dependent radius. The extension to Coulomb friction is briefly discussed at the end of Section~\ref{section:4:numres}.

\section{Model problems}\label{section:4:modpro}
Let $\mathcal{D} \subset\mathbb{R}^m, m \in \netoi:=\mathbb{N}\setminus\{0\}$, be the  parameter set. For all $\mu \in \mathcal{D}$, we consider an elastic body whose reference configuration is the parameter-dependent domain $\omegaa(\mu)\subset \mathbb{R}^{d}$, with $d \in \{2,3\}$. The boundary $\gammaa(\mu):=\partial\omegaa(\mu)$ is partitioned as $\gammaa(\mu) = \gammaad(\mu)\cup \gammaan(\mu) \cup \gammaac(\mu)$. The body is clamped at the boundary $\gammaad(\mu)$, free of traction at the boundary $\gammaan(\mu)$, and $\gammaac(\mu)$ denotes the potential contact boundary with a given rigid support. We denote by $\nug(\mu)$ the unit outward normal on $\gammaa(\mu)$ and by $\taug(\mu):= \big[\taug_{1}(\mu)\cdots \taug_{d-1}(\mu)\big] \in \Rd{d \times (d-1)}$ an orthonormal basis of the hyperplane orthogonal to $\nug(\mu)$ in $\Rd{d}$. For simplicity, we just write $\nug$ and $\taug$ whenever there is no ambiguity. The body in its reference configuration is located at some distance from a rigid support and we denote by $g(\mu)\in L^2(\gammaac(\mu);\mathbb{R}_{+})$ the corresponding gap function. An external load $\ellbm(\mu): \omegaa(\mu) \rightarrow \mathbb{R}^d$ is applied to the body, and we assume to be in the case of small deformations. For a generic $\mathbb{R}^{d}$-valued displacement field $\vbm$, the $\mathbb{R}^{d\times d}$-valued linearized strain tensor $\varepsilon(\vbm)$ and the $\mathbb{R}^{d\times d}$-valued stress tensor $\sigma(\vbm)$ are given by
\begin{align}
	\varepsilon(\vbm):=\frac{1}{2}\big(\nabla \vbm + \nabla \vbm^\top\big),\quad
	\sigma(\vbm):= \mathbb{C}\varepsilon(\vbm), \label{equa:4:stress}
\end{align}
with $\mathbb{C}$ the elastic coefficient tensor. At the boundary, we decompose the displacement field, $\vbm$, and the normal component of the stress tensor, $\sigma(\vbm)\nug$, in normal and tangential components as follows:
\begin{equation}
	\vbm =\vn\nug +  \taug\vt, \quad \sigma(\vbm)\nug=\sigman(\vbm)\nug + \taug\sigmat(\vbm),
\end{equation} 
with $\vn \in \mathbb{R}$, $\vt \in \Rd{d-1}$, $\sigman(\vbm) \in \mathbb{R}$ and $\sigmat(\vbm) \in \Rd{d-1}$.

The frictionless contact problem (also called Signorini problem) consists in finding the displacement field $\ubm(\mu):\omegaa(\mu) \rightarrow \mathbb{R}^d$ satisfying, for all $\mu \in \mathcal{D}$,
\begin{subequations}
	\label{equa:4:contact}
	\begin{align}
		-\mbox{div}(\sigma(\ubm(\mu)))  =  \ellbm(\mu),\quad&  \mbox{in } \omegaa(\mu), \label{equa:4:equi}\\ 
		\ubm(\mu) = \zerobm, \quad&  \mbox{on } \gammaad(\mu), \label{equa:4:diri}\\ 
		\sigma(\ubm(\mu))\nug = \zerobm, \quad&  \mbox{on } \gammaan(\mu), \label{equa:4:neum}\\
		\un(\mu) \leq g(\mu),\ \sigman(\ubm(\mu))\leq 0,\ \sigman(\ubm(\mu))(\un(\mu) - g(\mu)) =0, \quad& \mbox{on } \gammaac(\mu) \label{equa:4:signo},\\
		\sigmat(\ubm(\mu)) = \zerobm, \quad& \mbox{on } \gammaac(\mu).\label{equa:4:fricless}
	\end{align}
\end{subequations}

In the case of contact problems with friction, the condition~\eqref{equa:4:fricless} on $\gammaac(\mu)$ has to be replaced by a condition depending on the considered friction law~\cite{Kikuchi-1988,Mlika-2017}. Here, we focus on Tresca friction which leads to the following conditions:
\begin{equation}
\left\{
\begin{alignedat}{1}
\|\sigmat(\ubm(\mu))\| \leq s,&\quad \mbox{if }  \ut(\mu) =\zerobm,\\
\sigmat(\ubm(\mu)) = -s \frac{\ut(\mu)}{\|\ut(\mu)\|}, &\quad \mbox{otherwise},
\label{equa:4:tres}
\end{alignedat}
\right.
\end{equation}
where $\|\cdot\|$ denotes the Euclidean norm in $\mathbb{R}^{d-1}$ and  $s>0$ is a given threshold, taken to be constant for simplicity (units in Pa).
The model problem consisting of equations~\eqref{equa:4:equi}-\eqref{equa:4:diri}-\eqref{equa:4:neum}-\eqref{equa:4:signo} and \eqref{equa:4:tres} is called Tresca frictional contact problem.

We introduce the finite-dimensional space $\fespacebm{V}(\mu)$  and the admissible set $\enadm{K}(\mu)$ such that
\begin{subequations}
	\begin{align}
		\fespacebm{V}(\mu) &\subset\big\{\vbm \in  H^{1}(\omegaa(\mu);\mathbb{R}^{d}) \ |\  \vbm=\zerobm \mbox{ on } \gammaad(\mu)\big\}, \label{equa:4:vsp}\\
		\enadm{K}(\mu)&:=\big\{\vbm \in \fespacebm{V}(\mu)\ |\ \vn \leq g(\mu) \mbox{ on } \gammaac(\mu) \big\} \label{equa:4:kad}.
	\end{align}
\end{subequations}
We notice that $\enadm{K}(\mu)$ is a non-empty convex set. The space $\fespacebm{V}(\mu)$ is typically built as a finite element space associated with a fine mesh of $\omegaa(\mu)$.
The bilinear form $a(\mu;\cdot,\cdot): \fespacebm{V}(\mu) \times  \fespacebm{V}(\mu)  \rightarrow \mathbb{R}$  associated with the equilibrium equation~\eqref{equa:4:equi} in $\omegaa(\mu)$ is defined as
\begin{align}
	a(\mu;\ubm,\vbm):=\int_{\omegaa(\mu)} \sigma(\ubm):\varepsilon(\vbm)\,d\Omega(\mu),
\end{align}
and the linear form $f(\mu;\cdot): \fespacebm{V}(\mu)\rightarrow \mathbb{R}$ associated with the external load $\ellbm(\mu)$  as 
\begin{align}
	f(\mu;\vbm):=\int_{\omegaa(\mu)} \ellbm(\mu)\cdot \vbm\,  d\Omega(\mu).
\end{align}

The weak formulation of the Signorini contact problem \eqref{equa:4:contact} consists of solving the following variational inequality of the first kind: For all $\mu \in \mathcal{D}$,  find $\ubm(\mu) \in \enadm{K}(\mu)$ such that
\begin{align}
	a(\mu;\ubm(\mu),\vbm-\ubm(\mu)) \geq f(\mu;\vbm-\ubm(\mu)), \quad \forall \vbm \in \enadm{K}(\mu).
	\label{equa:4:invarsign}
\end{align}
For all $\mu \in \mathcal{D}$, Stampacchia's theorem~\cite{Stampacchia-1964} ensures that there is a unique solution to \eqref{equa:4:invarsign} which is also the unique solution to the following constrained minimization problem:  Find $\ubm(\mu) \in \enadm{K}(\mu)$ such that
\begin{align}
	\ubm(\mu) = \underset{\vbm \in \enadm{K}(\mu)}{\mbox{argmin}}\ \mathcal{J}(\mu;\vbm),
	\label{equa:4:VI} 
\end{align}
where the energy functional $\mathcal{J}(\mu;\cdot): \fespacebm{V}(\mu) \rightarrow \mathbb{R}$ is defined as follows:
\begin{align}
	\mathcal{J}(\mu;\vbm):= \frac{1}{2} a(\mu;\vbm,\vbm) - f(\mu;\vbm).
	\label{equa:4:energy}
\end{align}

In the case of Tresca friction, we need to consider the friction functional $\mathcal{F}$ such that
\begin{equation}
\mathcal{F}(\mu;\vbm):= 
\int_{\gammaac(\mu)}-s\|\vt\|\,d\gammaa(\mu).
\end{equation}
This leads to the following variational inequality: For all $\mu \in \mathcal{D}$, find $\ubm(\mu) \in \enadm{K}(\mu)$ such that
\begin{align}
a(\mu;\ubm(\mu),\vbm-\ubm(\mu)) + \mathcal{F}(\mu;\vbm) - \mathcal{F}(\mu;\ubm(\mu))\geq f(\mu;\vbm-\ubm(\mu)), \quad \forall \vbm \in \enadm{K}(\mu).
\label{equa:4:invarfr}
\end{align}

\section{Nitsche's method}\label{section:4:nitmet}
The main idea in the original Nitsche method~\cite{Nitsche-1971} is the enforcement of Dirichlet boundary conditions by means of a consistent penalty method. As shown in~\cite{Chouly-2013,Chouly-2015,Chouly-2017}, it is possible to generalize this idea to frictional contact problems. The main advantage of Nitsche's method is that the problem to be solved is unconstrained. Hence, in contrast to the mixed formulation, one does not need any additional unknowns such as Lagrange multipliers. Moreover, a higher effectivity of the RBM is expected since it is well-known that the dual basis is particularly hard to compress. The price to be paid, though, is that the constraint is not exactly enforced. 

To derive Nitsche's method, we are going to assume that all the considered functions are smooth enough so that the associated normal stress tensor can be defined pointwise at the boundary. Recall that, for all $\mu \in \mathcal{D}$, the HF finite-dimensional space $\fespacebm{V}(\mu)$ in~\eqref{equa:4:vsp} results from a finite element discretization of the Hilbert space $H^{1}(\Omega(\mu);\Rd{d})$. Hence, the above assumption is indeed met.
\subsection{Frictionless case}\label{subsection:4:nofric}
Following~\cite{Curnier-1988}, the starting observation is that the Signorini conditions~\eqref{equa:4:signo} can be equivalently reformulated as follows:
\begin{align}
	\sigman(\ubm(\mu)) = \projconv{\sigman(\ubm(\mu)) - \gamma \big(\un(\mu) - g(\mu)\big)}{-},
	\label{equa:4:nitsrefor}
\end{align}
where $\projconv{z}{-}:=\mbox{min}\big(z,0\big)$ denotes the negative part of a generic real number $z$ and where $\gamma>0$ is a user-defined parameter (taken to be constant for simplicity). In practice, the parameter $\gamma$ should be chosen large enough (see Section~\ref{section:4:numres} for further discussion). 

For all $\mu \in \mathcal{D}$, one introduces  the  energy functional $\JN(\mu;\cdot): \fespacebm{V}(\mu) \rightarrow \mathbb{R}$ such that
\begin{equation}
	\begin{alignedat}{1}
		\JN(\mu; \vbm):=&\  \mathcal{J}(\mu;\vbm) - \frac{1}{2} \int_{\gammaac(\mu)} \frac{1}{\gamma}| \sigman(\vbm)|^2 d\gammaa(\mu)\\
		&+\frac{1}{2} \int_{\gammaac(\mu)} \frac{1}{\gamma} \projconv{\sigman(\vbm) -\gamma \big(\vn - g(\mu)\big)}{-}^2 d\gammaa(\mu),
	\end{alignedat}
\end{equation}
recalling that the energy functional $\mathcal{J}$ is defined in~\eqref{equa:4:energy}.
Nitsche's method consists in finding $\ubm(\mu) \in \fespacebm{V}(\mu)$ solution to the following unconstrained minimization problem: Find $\ubm(\mu) \in \fespacebm{V}(\mu)$ such that
\begin{align}
	\ubm(\mu) = \underset{\vbm \in \fespacebm{V}(\mu)}{\mbox{argmin}}\ \JN(\mu;\vbm).
	\label{equa:4:VIN} 
\end{align}
The first-order optimality condition associated with \eqref{equa:4:VIN} reads 
\begin{equation}
	\Agamn(\mu;\ubm(\mu),\vbm) + \int_{\gammaac(\mu)} \frac{1}{\gamma} \projconv{\Pgamn(\mu;\ubm(\mu))}{-} \Pgamnz(\mu;\vbm)\, d\gammaa(\mu) = f(\mu;\vbm),\quad \forall \vbm \in \fespacebm{V}(\mu),
	\label{equa:4:NWOF}
\end{equation}
with the bilinear form $\Agamn(\mu;\cdot,\cdot) : \fespacebm{V}(\mu) \times \fespacebm{V}(\mu) \rightarrow \mathbb{R}$ defined as
\begin{equation}
	\Agamn(\mu;\ubm,\vbm):=a(\mu;\ubm,\vbm) -  \int_{\gammaac(\mu)} \frac{1}{\gamma} \sigman(\ubm)\sigman(\vbm)\,d\Gamma(\mu),\quad \forall \ubm, \vbm \in \fespacebm{V}(\mu),
\end{equation}
and the operators $\Pgamn(\mu;\cdot),\ \Pgamnz(\mu;\cdot): \fespacebm{V}(\mu)\rightarrow L^2(\gammaac(\mu))$  defined as 
\begin{subequations}
	\begin{alignat}{2}
		&\Pgamn(\mu;\vbm):=\sigman(\vbm) -\gamma(\vn-g(\mu)),\\
		&\Pgamnz(\mu;\vbm):=\sigman(\vbm) -\gamma \vn.
	\end{alignat}
\end{subequations}
With this notation, we can rewrite the  energy functional $\JN(\mu;\cdot)$ as
\begin{equation}
	\begin{alignedat}{1}
		\JN(\mu; \vbm):=&\  \mathcal{J}(\mu;\vbm) - \frac{1}{2} \int_{\gammaac(\mu)} \frac{1}{\gamma} |\sigman(\vbm)|^2 d\gammaa(\mu) +
		\frac{1}{2} \int_{\gammaac(\mu)} \frac{1}{\gamma} \projconv{\Pgamn(\vbm)}{-}^2 \,d\gammaa(\mu).
	\end{alignedat}
\end{equation}

The problem~\eqref{equa:4:NWOF} is nonlinear. To solve it, we use an iterative method. Given $\ubm_{0}(\mu) \in \fespacebm{V}(\mu)$ and a tolerance $ \prpu \in \rdp$, we look, for all $k\geq 0$, for the solution $\ubm_{k+1}(\mu)$  in the form $\ubm_{k+1}(\mu) = \ubm_{k}(\mu) + \deltaubm_{k}(\mu)$. Ideally, we seek $\deltaubm_{k}(\mu) \in \fespacebm{V}(\mu)$ such that 
\begin{equation}
	\left\{
	\begin{alignedat}{2}
		&\Agamn(\mu;\ubm_{k}(\mu) + \deltaubm_{k}(\mu),\vbm) \\
		&+ \int_{\gammaac(\mu)} \frac{1}{\gamma} \projconv{\Pgamn(\mu;\ubm_{k}(\mu) + \deltaubm_{k}(\mu))}{-} \Pgamnz(\mu;\vbm)\, d\gammaa(\mu) 
		= f(\mu;\vbm),\quad \forall \vbm \in \fespacebm{V}(\mu).
	\end{alignedat}
	\right.
	\label{equa:4:itNWOF}
\end{equation}
However, since the problem~\eqref{equa:4:itNWOF} is nonlinear, it is expensive to solve it directly. Therefore, we approximate the solution $\deltaubm_{k}(\mu)$ by linearizing the problem. In order to do this, we observe that $\Pgamn(\mu;\ubm_{k}(\mu) + \deltaubm_{k}(\mu)) = \Pgamn(\mu;\ubm_{k}(\mu))+ \Pgamnz(\mu; \deltaubm_{k}(\mu))$ and  approximate the term $\projconv{\Pgamn(\mu;\ubm_{k}(\mu) + \deltaubm_{k}(\mu))}{-}$  as follows: 
\begin{align}
	\projconv{\Pgamn(\mu;\ubm_{k}(\mu) + \deltaubm_{k}(\mu))}{-} 
	\approx \projconv{\Pgamn(\mu;  \ubm_{k}(\mu))}{-} + H(-\Pgamn(\mu; \ubm_{k}(\mu))) \Pgamnz(\mu; \deltaubm_{k}(\mu)),
\end{align}
where $H(\cdot)$ denotes the Heaviside function. Using this approximation in~\eqref{equa:4:itNWOF} (for simplicity, we keep the same notation for the unknown $\deltaubm_{k}(\mu)$), we consider the following sequence of problems: For all $k\geq 0$, find $\deltaubm_{k}(\mu) \in \fespacebm{V}(\mu)$ such that
\begin{equation}
	\Agamn(\mu; \deltaubm_{k}(\mu),\vbm) + \Bgamn(\mu;\ubm_{k}(\mu); \deltaubm_{k}(\mu),\vbm) = -\Resn(\mu;\ubm_{k}(\mu);\vbm),\quad \forall \vbm \in \fespacebm{V}(\mu),
	\label{equa:4:itNWOF2}
\end{equation}
where 
\begin{subequations}
	\begin{alignat}{2}
		\Bgamn(\mu;\wbm;\ubm,\vbm) &:= \int_{\gammaac(\mu)}\frac{1}{\gamma} H(-\Pgamn(\mu; \wbm)) \Pgamnz(\mu; \ubm)\Pgamnz(\mu; \vbm)\,d\gammaa(\mu),\label{equa:4:opb}\\
		\Resn(\mu;\wbm;\vbm)&:=  \Agamn(\mu; \wbm,\vbm) + \Rhon(\mu;\wbm;\vbm) - f(\mu;\vbm),\label{equa:4:residu}\\
	\end{alignat}
\end{subequations}
with 
\begin{equation}
	\Rhon(\mu;\wbm,\vbm):= \int_{\gammaac(\mu)}\frac{1}{\gamma} \projconv{\Pgamn(\mu; \wbm)}{-}\Pgamnz(\mu; \vbm)\,d\gammaa(\mu).\label{equa:4:thetan}
\end{equation}
We iterate on $k$ until the following convergence criterion is reached: 
\begin{equation}
	\frac{\norm{\ubm_{k+1}(\mu) - \ubm_{k}(\mu)}{\fespacebm{V}(\mu)}}{ \norm{\ubm_{k+1}(\mu)}{\fespacebm{V}(\mu)}} \leq \prpu.
	\label{equa:4:cvcrit}
\end{equation}
In what follows, we denote by $\kc(\mu) \in \netoi$ the number of iterations required to solve~\eqref{equa:4:itNWOF2}, for all $\mu \in \mathcal{D}$.
We denote the converged solution to the sequence of problems~\eqref{equa:4:itNWOF2} as $\ubm_{\rm{cv}}(\mu) \in \fespacebm{V}(\mu)$.

To introduce the algebraic formulation, we assume that for all $\mu \in \mathcal{D}$, the HF (finite-dimensional) space $\fespacebm{V}(\mu)$ is such that 
\begin{align}
	\fespacebm{V}(\mu)&:= \bspan\Iparr{\setcalmuhf{\varphibm}{n}{1}{\mu}}.
\end{align}
Notice that the dimension $\dimP$ of $\fespacebm{V}(\mu)$ is parameter-independent. We illustrate in Section~\ref{subsection:4:geom} how to accomplish this property. 
Furthermore, we adopt the following decompositions for the sequence of solutions to~\eqref{equa:4:itNWOF2}:
\begin{subequations}
	\label{equa:4:hfreconstr}
	\begin{alignat}{4}
		\ubm_{k}(\mu)&:=\Sum{n}{1}{\dimP}\ubmv{n}_{k}(\mu)\varphibm_{n}(\mu),\quad & \Ubm_{k}(\mu)&:=\Vect{\ubmv{n}_{k}}{\mu}{n}{1}{\dimP} \in \Rd{\dimP},\\
		\deltaubm_{k}(\mu)&:=\Sum{n}{1}{\dimP}\deltaubmv{n}_{k}(\mu)\varphibm_{n}(\mu),\quad & \Deltaubm_{k}(\mu)&:=\Vect{\deltaubmv{n}_{k}}{\mu}{n}{1}{\dimP} \in \Rd{\dimP}.
	\end{alignat}
\end{subequations}
The algebraic formulation of the sequence of problems~\eqref{equa:4:itNWOF2} then reads as follows: For all $k\geq 0$, find $\Deltaubmka{k}(\mu) \in \mathbb{R}^{\dimP}$ such that
\begin{equation}
	\Agamnm(\mu)  \Deltaubmka{k}(\mu) + \Bgamnm(\mu,\ubmka{k}(\mu))  \Deltaubmka{k}(\mu) = - \Resnm(\mu,\ubmka{k}(\mu)),
	\label{equa:4:itNWOF2Alg}
\end{equation}
where for all $n,m \in \indset{1}{\dimP}$,
\begin{subequations}
	\begin{alignat}{2}
		(\Agamnm(\mu))_{mn}&:=\Agamn(\mu;\varphi_{n}(\mu),\varphi_{m}(\mu)),\\
		(\Bgamnm(\mu,\ubmka{k}(\mu)))_{mn}&:=\Bgamn(\mu;\ubm_{k}(\mu);\varphi_{n}(\mu),\varphi_{m}(\mu)),  \\
		(\Resnm(\mu,\ubmka{k}(\mu)))_{m}&:=\Resn(\mu;\ubm_{k}(\mu);\varphi_{m}(\mu)).
	\end{alignat}
\end{subequations}
The convergence criterion for~\eqref{equa:4:itNWOF2Alg} is still~\eqref{equa:4:cvcrit} using the reconstructed functions (see~\eqref{equa:4:hfreconstr}). We denote the converged solution to the sequence of problems~\eqref{equa:4:itNWOF2Alg} as $\Ubm_{\rm{cv}}(\mu) \in \Rd{\dimP}$. Furthermore, let us denote by $\Rhonm(\mu,\wbm)$ (resp. $F(\mu)$) the algebraic representation of the operator $\Rhon(\mu;\wbm;\cdot)$ (resp. $f(\mu;\cdot)$) such that for all $ m \in \indset{1}{\dimP}$,
\begin{equation}
	(\Rhonm(\mu,\wbm))_{m}:= \Rhon(\mu;\wbm;\varphi_{m}(\mu)), \qquad (F(\mu))_{m}:= f(\mu;\varphi_{m}(\mu)).
\end{equation}
With this notation, we have the following decomposition:
\begin{equation}
	\Resnm(\mu,\ubmka{k}(\mu)):= \Agamnm(\mu)\Ubmka{k}(\mu) + \Rhonm(\mu,\ubmka{k}(\mu))-F(\mu).
\end{equation}

\subsection{Friction case}\label{subsection:4:fric}
For the Tresca frictional contact problem, we use the following reformulation of the friction conditions given in~\cite{Chouly-2017}:
\begin{align}
	\sigmat(\ubm(\mu)) = \projconv{\sigmat(\ubm(\mu)) - \gamma \ut(\mu)}{s},
\end{align}
where, for a  positive real number $r$, the notation
\begin{equation}
	\projconv{\xbm}{r}:= 
	\begin{cases}
		\xbm,&\quad \mbox{if } \|\xbm\|\leq r,\\
		r\frac{\xbm}{\| \xbm\|}, &\quad \mbox{otherwise},
	\end{cases}
\end{equation}
defines the projection of $\xbm \in \mathbb{R}^{d-1}$ onto the ball centered at the origin and of radius $r$. Let us introduce the operator $\Pgamt(\mu;\cdot):\fespacebm{V}(\mu) \rightarrow L^2(\gammaac(\mu);\mathbb{R}^{d-1})$ such that 
\begin{align}
	\Pgamt(\mu;\vbm):=\sigmat(\vbm) -\gamma \vt.
\end{align}
We start from the  energy functional $\JNF(\mu;\cdot): \fespacebm{V}(\mu) \rightarrow \mathbb{R}$ defined as
\begin{equation}
	\begin{alignedat}{2}
		\JNF(\mu;\vbm):=&~\JN(\mu;\vbm) - \frac{1}{2} \int_{\gammaac(\mu)} \frac{1}{\gamma}\| \sigmat(\vbm)\|^2\, d\gammaa(\mu) \\
		&+ \frac{1}{2} \int_{\gammaac(\mu)} \frac{1}{\gamma} \| \Pgamt(\vbm)\|^{2}\, d\gammaa(\mu)
		-\frac{1}{2} \int_{\gammaac(\mu)} \frac{1}{\gamma}\| \Pgamt(\vbm) -\projconv{\Pgamt(\vbm)}{s} \|^{2}\, d\gammaa(\mu).
	\end{alignedat}
	\label{equa:4:NWFE}
\end{equation}
Nitsche's method consists in solving the following unconstrained minimization problem: Find $\ubm(\mu) \in \fespacebm{V}(\mu)$ such that
\begin{align}
	\ubm(\mu) = \underset{\vbm \in \fespacebm{V}(\mu)}{\mbox{argmin}}\ \JNF(\mu;\vbm).
	\label{equa:4:VINWF} 
\end{align}
Let us consider the functional $J:\mathbb{R}^{d-1} \rightarrow \mathbb{R}$ defined as follows:
\begin{equation}
	J(\xbm) = \frac{1}{2}\|\xbm - \projconv{\xbm}{s} \|^{2}.
\end{equation}
One easily verifies that this functional is Gâteaux-differentiable and  that its differential is given by $DJ(\xbm)= \xbm -\projconv{\xbm}{s}$ for all $\xbm \in \mathbb{R}^{d-1}$.
Therefore, using the same approach as for the frictionless contact problem, we obtain the following Nitsche's formulation for the Tresca frictional contact problem:
Find $\ubm(\mu) \in \fespacebm{V}(\mu)$ such that
\begin{equation}
	\left\{
	\begin{alignedat}{2}
		&	\Agam(\mu;\ubm(\mu),\vbm) 
		+  \int_{\gammaac(\mu)} \frac{1}{\gamma} \projconv{\Pgamn(\mu;\ubm(\mu))}{-} \Pgamnz(\mu;\vbm)\, d\gammaa(\mu)\\
		& + \int_{\gammaac(\mu)} \frac{1}{\gamma} \projconv{\Pgamt(\mu;\ubm(\mu))}{s} \Pgamt(\mu;\vbm)\, d\gammaa(\mu) 
		=f(\mu;\vbm), \quad \forall \vbm \in \fespacebm{V}(\mu),
	\end{alignedat}
	\right.
	\label{equa:4:NWF}
\end{equation}
where the bilinear form $\Agam(\mu;\cdot,\cdot): \fespacebm{V}(\mu) \times \fespacebm{V}(\mu)\rightarrow \mathbb{R}$ is defined as 
\begin{equation}
	\Agam(\mu;\ubm,\vbm):=a(\mu;\ubm,\vbm) -  \int_{\gammaac(\mu)} \frac{1}{\gamma} \sigma(\ubm)\nug\cdot\sigma(\vbm)\nug \,d\Gamma(\mu).
\end{equation}
The linearization of the problem~\eqref{equa:4:NWF} leads to the following sequence of problems: For all $k\geq 0$, find $\deltaubm_{k}(\mu) \in \fespacebm{V}(\mu)$ such that 
\begin{equation}
	\Agam(\mu; \deltaubm_{k}(\mu),\vbm) + \Bgam(\mu;\ubm_{k}(\mu); \deltaubm_{k}(\mu),\vbm) = -\Res(\mu;\ubm_{k}(\mu);\vbm),\quad \forall \vbm \in \fespacebm{V}(\mu),
	\label{equa:4:itNWF}
\end{equation}
with
\begin{subequations}
	\begin{alignat}{2}
		\Bgam(\mu;\wbm;\ubm,\vbm)&:= \Bgamn(\mu;\wbm;\ubm,\vbm) +\Bgamt(\mu;\wbm;\ubm,\vbm),\label{equa:4:opbf} \\
		\Res(\mu;\wbm;\vbm)&:= \Agam(\mu; \wbm,\vbm) + \Rhon(\mu;\wbm,\vbm) + \Rhot(\mu;\wbm,\vbm)- f(\mu;\vbm),\label{equa:4:residuf}
	\end{alignat}
\end{subequations}
where
\begin{subequations}
\begin{alignat}{2}
\Bgamt(\mu;\wbm;\ubm,\vbm)&:=\int_{\gammaac(\mu)}\frac{1}{\gamma} \Big(\Gdif{s}(\Pgamt(\mu; \wbm)) \cdot\Pgamt(\mu; \ubm)\Big)\Pgamt(\mu; \vbm)\,d\gammaa(\mu),\\
	\Rhot(\mu;\wbm,\vbm)&:= \int_{\gammaac(\mu)}\frac{1}{\gamma} \projconv{\Pgamt(\mu; \wbm)}{s}\Pgamt(\mu; \vbm)\,d\gammaa(\mu),
\end{alignat}
\end{subequations}
and $\Gdif{s}(\cdot)$, the differential of $\projconv{\cdot}{s}$, is given by
\begin{equation}
	\Gdif{s}(\xbm):= 
	\begin{cases}
		\mathbb{I}_{d-1},&\quad \mbox{if } \|\xbm\|\leq s,\\
		\frac{s}{\| \xbm\|}\big(\mathbb{I}_{d-1} - \frac{\xbm \otimes \xbm}{\norm{\xbm}{}^2}\big), &\quad \mbox{otherwise},
	\end{cases}
\end{equation}
with $\mathbb{I}_{d-1}$ the identity matrix of order $d-1$. The convergence criterion for~\eqref{equa:4:itNWF} is still~\eqref{equa:4:cvcrit}. We denote the converged solution to the sequence of problems~\eqref{equa:4:itNWF} as $\ubm_{\rm{cv}}(\mu) \in \fespacebm{V}(\mu)$.

The algebraic formulation of the sequence of problems~\eqref{equa:4:itNWF} reads as follows: For all $k\geq 0$, find $\Deltaubmka{k}(\mu) \in \mathbb{R}^{\dimP}$ such that
\begin{equation}
	\Agamm(\mu)  \Deltaubmka{k}(\mu) + \Bgamm(\mu,\ubmka{k}(\mu))  \Deltaubmka{k}(\mu) = - \Resm(\mu,\ubmka{k}(\mu)),
	\label{equa:4:itNWF2Alg}
\end{equation}
where for all $n,m\in \indset{1}{\dimP}$,
\begin{subequations}
	\begin{alignat}{2}
		(\Agamm(\mu))_{mn}&:=\Agam(\mu;\varphi_{n}(\mu),\varphi_{m}(\mu)),\\
		(\Bgamm(\mu,\ubmka{k}(\mu) ))_{mn}&:=\Bgam(\mu;\ubm_{k}(\mu);\varphi_{n}(\mu),\varphi_{m}(\mu)),  \\
		(\Resm(\mu,\ubmka{k}(\mu) ))_{m}&:=\Res(\mu;\ubm_{k}(\mu);\varphi_{m}(\mu)).
	\end{alignat}
\end{subequations}
The convergence criterion for~\eqref{equa:4:itNWF2Alg} is still~\eqref{equa:4:cvcrit} using the reconstructed functions. We denote the converged solution to the sequence of problems~\eqref{equa:4:itNWF2Alg} as $\Ubm_{\rm{cv}}(\mu) \in \Rd{\dimP}$. Furthermore, let us denote by $\Bgamtm(\mu,\wbm)$ (resp. $\Rhotm(\mu,\wbm)$) the algebraic representation of the operator $\Bgamt(\mu;\wbm;\cdot,\cdot)$ (resp. $\Rhot(\mu;\wbm;\cdot)$) such that for all $n,m \in \indset{1}{\dimP}$,
\begin{equation}
(\Bgamtm(\mu,\wbm))_{mn}:=\Bgamt(\mu;\wbm;\varphi_{n}(\mu),\varphi_{m}(\mu)),\qquad (\Rhotm(\mu,\wbm))_{m}:= \Rhot(\mu;\wbm;\varphi_{m}(\mu)).
\end{equation}
With this notation, we have the following decompositions:
\begin{subequations}
	\begin{alignat}{2}
		\Bgamm(\mu,\ubmka{k}(\mu))&:= \Bgamnm(\mu,\ubmka{k}(\mu)) + \Bgamtm(\mu,\ubmka{k}(\mu)),\\
		\Resm(\mu,\ubmka{k}(\mu))&:= \Agamm(\mu)\Ubmka{k}(\mu) + \Rhonm(\mu,\ubmka{k}(\mu))+\Rhotm(\mu,\ubmka{k}(\mu))-F(\mu).
	\end{alignat}
\end{subequations}

\section{Reduced-basis formulation}\label{section:4:rbm}
In this section, we derive the RBM. We describe the problem in detail for the frictionless contact problem, and briefly highlight the (simple) adaptations needed to account for friction.
\subsection{Geometric mapping}\label{subsection:4:geom}
We recall that the dimension $\dimP$ of the HF (finite-dimensional) space $\fespacebm{V}(\mu)$ is parameter-independent. This is important since, in order to compress the space generated by the snapshots, it is necessary that all the snapshots live in the same space. For this purpose, since the geometry is parameter-dependent, we use a parameter-independent reference domain $\omegaahat$. We assume that for all $\mu \in \mathcal{D}$, there is a smooth invertible geometric mapping $\oph(\mu): \omegaahat \rightarrow \omegaa(\mu)$. We denote by $\gammaahat :=\partial\omegaahat$ the boundary of $\omegaahat$ and we assume that it can be partitioned as $\gammaahat = \gammaadhat \cup \gammaanhat \cup \gammaachat$  in such a way that, for all $\mu \in \mathcal{D}$,
\begin{align}
	\gammaadhat:=\ophd^{-1}(\mu)(\gammaad(\mu)), \quad   \gammaanhat:=\ophn^{-1}(\mu)(\gammaan(\mu)),\quad  \gammaachat:=\ophc^{-1}(\mu)(\gammaac(\mu)),
\end{align}  
with $\ophd(\mu):=\oph(\mu)_{|\gammaadhat}$,  $\ophn(\mu):=\oph(\mu)_{|\gammaanhat}$ and $\ophc(\mu):=\oph(\mu)_{|\gammaachat}$. Then, the mesh of $\omegaa(\mu)$ is generated by generating a mesh of $\omegaahat$ matching the partition of the boundary $\gammaahat$ and applying the mapping $\oph(\mu)$ to the mesh of the reference domain $\omegaahat$. The finite element basis functions are generated from the reference basis functions by using a pullback.

\subsection{Naive approach}\label{subsection:4:NA}
The goal of this section is to present the naive reduced model resulting from the application of a plain RBM to the contact problem formulated with Nitsche's method and highlight the computational inefficiency of such a formulation. This problem will be circumvented in the next section eventually leading to a computationally effective RBM. The major difficulty comes from the nonlinearity of Nitsche's formulation.

To build the reduced basis (RB), the starting point is to compute (in the offline phase) a family $\setbis{\Ubmka{\rm{cv}}(\mup{p})}{p}{1}$ $\subset  \Rd{\dimP}$ of HF solutions to the frictionless contact problem~\eqref{equa:4:itNWOF2Alg} by using a training subset $\discrd{D}{train}:=\set{\mu}{p}{1} \subset \mathcal{D}$ of cardinality $P \in \netoi$. Using the Proper Orthogonal Decomposition (POD)~\cite{Haasdonk-2013,Kunisch-2001} based on the canonical inner product of $H^{1}(\omegaahat;\Rd{d})$ and the geometric mapping $\oph(\mu)$, one can construct an orthonormal family $\setbbb{\Xibm}{n}{1} \subset \Rd{\dimP}$ of $N\in \netoi$  ($N \leq P$) vectors. Let us denote by $\rfespacebm{V}{N}$ the reduced space generated by the family $\setbbbbis{\Xibm}{n}{1}$, \ie,
\begin{equation}
	\rfespacebm{V}{N}:= \bspan\Iparr{\setbbbbis{\Xibm}{n}{1} } \subset \Rd{\dimP}.
\end{equation}

In algebraic form, the RB formulation of the sequence of HF problems~\eqref{equa:4:itNWOF2Alg} reads as follows: 
For all $k\geq 0$, find $\Deltaubmredka{k}(\mu) \in \mathbb{R}^{\dimPred}$ such that
\begin{equation}
	\AgamnM(\mu) \Deltaubmredka{k}(\mu) + \BgamnM(\mu,k) \Deltaubmredka{k}(\mu) = - \ResnM(\mu,k),
	\label{equa:4:itNWOF2refred}
\end{equation}
where 
\begin{subequations}
	\begin{alignat}{4}
		\AgamnM(\mu)&:= \Z^\top\Agamnm(\mu)\Z \in \Rd{\dimPred\times \dimPred},\label{equa:4:RBmatrA}\\
		\BgamnM(\mu,k)&:= \Z^\top\Bgamnm(\mu,\ubmredka{k}(\mu))\Z \in \Rd{\dimPred\times \dimPred}, \quad
		\ResnM(\mu,k)&:=\Z^\top\Resnm(\mu,\ubmredka{k}(\mu)) \in \Rd{\dimPred},\label{equa:4:RBmatrBR}
	\end{alignat}
	\label{equa:4:RBmatr}
\end{subequations}
with $\Z:=\matr{\Xibm_{1}}{\Xibm_{\dimPred}} \in \Rd{\dimP\times \dimPred}$.
The convergence criterion for~\eqref{equa:4:itNWOF2refred} is still~\eqref{equa:4:cvcrit} using the reconstructed functions (see Remark~\ref{remark:4:funcreconstr} for more details).

At this stage, the RBM consists of the following two stages:
\begin{itemize}
	\item[---] Offline stage
	\begin{enumerate}
		\item Select a training subset $\discrd{D}{train}:=\set{\mu}{p}{1} \subset \mathcal{D}$.
		\item Compute the snapshots $\setbis{\Ubmka{\rm{cv}}(\mup{p})}{p}{1} \subset \Rd{\dimP}$ by solving~\eqref{equa:4:itNWOF2Alg}.
		\item Compute the reduced space $\rfespacebm{V}{N}$ by using POD on snapshots.
	\end{enumerate}
	\item[---] Online stage:  For any $\mu \in \mathcal{D}\setminus \discrd{D}{train}$,
	\begin{enumerate}
		\item Compute $\AgamnM(\mu)$ using~\eqref{equa:4:RBmatrA}.
		\item Loop on $k\geq 0$,
		\begin{enumerate}
			\item \label{item:4:RBalgo2} Compute $\BgamnM(\mu,k)$ and $\ResnM(\mu,k)$ using~\eqref{equa:4:RBmatrBR}.
			\item Solve~\eqref{equa:4:itNWOF2refred}. 
			\item Check convergence; if not, set $k \leftarrow k + 1$ and go back to~(\ref{item:4:RBalgo2}).
		\end{enumerate}
	\end{enumerate}
\end{itemize}
It is crucial to derive a reduced problem that is independent of the high-fidelity dimension $\dimP$ in order to obtain an inexpensive online stage. This condition is not yet satisfied with the current formalism. The main issue is the manipulation of large-dimensional arrays in~\eqref{equa:4:itNWOF2refred}. We propose in Section~\ref{subsection:4:RBMEIM} a procedure to overcome this issue in order to construct a computationally efficient RBM.

\begin{remark}[Reconstructed functions]\label{remark:4:funcreconstr}
	Let us introduce the following reconstructed functions:
	\begin{equation}
		\begin{alignedat}{3}
			\xibm_{n}(\mu)&:= \Sum{i}{1}{\dimP} \Xibm_{n}^{i}\varphibm_{i}(\mu) \in \fespacebm{V}(\mu),\quad &\forall n \in \indset{1}{\dimPred} .
		\end{alignedat}
	\end{equation}
	With this notation, solving the RB problem~\eqref{equa:4:itNWOF2refred} in algebraic form leads to the following reconstructed solutions:
	\begin{subequations}
		\begin{alignat}{4}
			\ubmredka{k}(\mu)&:=\Sum{n}{1}{N}\ubmredvka{n}{k}(\mu)\xibm_{n}(\mu) \in \fespacebm{V}(\mu),\quad & \Ubmredka{k}(\mu)&:=\Vect{\ubmredvka{n}{k}}{\mu}{n}{1}{N} \in \mathbb{R}^{N}, \\
			\deltaubmredka{k}(\mu)&:=\Sum{n}{1}{N}\deltaubmredvka{n}{k}(\mu)\xibm_{n}(\mu)\in \fespacebm{V}(\mu),\quad & \Deltaubmredka{k}(\mu)&:=\Vect{\deltaubmredvka{n}{k}}{\mu}{n}{1}{N} \in \mathbb{R}^{N}.
		\end{alignat}
	\end{subequations}
\end{remark}

\subsection{Computationally efficient approach} \label{subsection:4:RBMEIM}
In this section, we describe the strategy to avoid the manipulation of large-dimensional arrays in the problem~\eqref{equa:4:itNWOF2refred}. The idea consists in introducing appropriate affine parametric decompositions of the parameter-dependent and ``parameter/iteration''-dependent operators involved in the problem by using the Empirical Interpolation Method (EIM)~\cite{Barrault-2004,Maday-2009}. 

Specifically, our goal is to separate the dependence on $\mu$ and $k$ from the dependence on the indices in the large-dimensional arrays {$\Agamnm(\mu)$, $\Bgamnm(\mu,\ubmredka{k}(\mu))$ and $\Resnm(\mu,\ubmredka{k}(\mu))$.} This operation is performed during the offline stage. For this purpose, using the EIM, we obtain the following approximations:
\begin{subequations}
	\label{equa:4:EIMapp}
	\begin{alignat}{5}
		\Agamnm(\mu)&\approx&\eimop{{a^{\nug}}}(\mu)&:=\Sum{s}{1}{\eimcar{{a^{\nug}}}} \Alphas{s}{{a^{\nug}}}(\mu)\Agamns,&\quad& \Agamns \in \mathbb{R}^{\dimP\times \dimP}, \, \Alphas{s}{{a^{\nug}}}(\mu)\in \mathbb{R},\label{equa:4:EIMappA}\\
		\Bgamnm(\mu,\ubmka{k}(\mu))&\approx{}&\eimop{{b^{\nug}}}(\mu,k)&:=\Sum{s}{1}{\eimcar{{b^{\nug}}}} \Alphas{s}{{b^{\nug}}}(\mu,k)\Bgamns,&\quad& \Bgamns \in \mathbb{R}^{\dimP\times \dimP},\, \Alphas{s}{{b^{\nug}}}(\mu,k)\in \mathbb{R},\label{equa:4:EIMappB}\\
		\Rhonm(\mu,\ubmka{k}(\mu))&\approx{}&{\eimop{{\theta^{\nug}}}(\mu,k)}&:=\Sum{s}{1}{\eimcar{{\theta^{\nug}}}} \Alphas{s}{{\theta^{\nug}}}(\mu,k)\Rhons,&\quad& \Rhons \in \mathbb{R}^{\dimP}, \, \Alphas{s}{{\theta^{\nug}}}(\mu,k)\in \mathbb{R},\label{equa:4:EIMappTheta}\\
		\F(\mu)&\approx&\eimop{f}(\mu)&:=\Sum{s}{1}{\eimcar{f}} \Alphas{s}{f}(\mu)\Fs,&\quad& \Fs \in \mathbb{R}^{\dimP}, \, \Alphas{s}{f}(\mu)\in \mathbb{R},\label{equa:4:EIMappF}
	\end{alignat}
\end{subequations} 
where the large-dimensional arrays $\Agamns$, $\Bgamns$, $\Rhons$ and $\Fs$ are now independent of the parameter $\mu$ and the iteration counter $k$, whereas the functions $\Alphas{s}{{a^{\nug}}}$ and ${\Alphas{s}{f}}$ (resp.,~$\Alphas{s}{{b^{\nug}}}$ and $\Alphas{s}{{\theta^{\nug}}}$) only depend on $\mu$ (resp.,~$(\mu,k)$). {We obtain the following approximation of the residual:
	\begin{equation}		\Resnm(\mu,\ubmka{k}(\mu))\approx{}\eimop{{r^{\nug}}}(\mu,k):=\eimop{{a^{\nug}}}(\mu)\Ubmka{k}(\mu)+\eimop{{\theta^{\nug}}}(\mu,k)-\eimop{f}(\mu). \label{equa:4:EIMappR}
\end{equation}}
To build the large-dimensional arrays $\Agamns$, $\Bgamns$,  $\Rhons$ and $\Fs$ , we respectively use the training sets $\eimtr{{a^{\nug}}}$, $\eimtr{{b^{\nug}}}$,  {$\eimtr{{\theta^{\nug}}}$} and {$\eimtr{f}$} defined as follows:
\begin{align}
	\eimtr{{a^{\nug}}}={\eimtr{f}}:=\discrd{D}{train}, \quad 	\eimtr{{b^{\nug}}}={\eimtr{{\theta^{\nug}}}}:= \big\{(\mu,k)~|~\mu \in \discrd{D}{train},~k \in  \indset{1}{\kc(\mu)}\big\}.
	\label{equa:4:eimtrain}
\end{align}
Notice that a different training set (possibly richer) than $\discrd{D}{train}$ can be used instead.
We introduce the index subsets $\{(\is{{a^{\nug}}},\js{{a^{\nug}}})\}_{s \in \indset{1}{\eimcar{{a^{\nug}}}}}, \{(\is{{b^{\nug}}},\js{{b^{\nug}}})\}_{s \in \indset{1}{\eimcar{{b^{\nug}}}}} \subset \indset{1}{\dimP} \times \indset{1}{\dimP}$, {$\{\js{{\theta^{\nug}}}\}_{s \in \indset{1}{\eimcar{{\theta^{\nug}}}}} \subset \indset{1}{\dimP}$}  and {$\{\js{f}\}_{s \in \indset{1}{\eimcar{f}}} \subset \indset{1}{\dimP}$} of cardinality $\eimcar{{a^{\nug}}}$, $\eimcar{{b^{\nug}}}$, {$\eimcar{{\theta^{\nug}}}$}, and {$\eimcar{f}$} respectively, corresponding to the indices selected by the EIM for the approximation of  $\Agamnm$, $\Bgamnm$, {$\Rhonm$} and {$\F$}, respectively. 
Then, the functions $\Alphas{s}{{a^{\nug}}}$, $\Alphas{s}{{b^{\nug}}}$,  {$\Alphas{s}{{\theta^{\nug}}}$} and {$\Alphas{s}{f}$} are such that 
\begin{subequations}
	\label{equa:4:eimcoef}
	\begin{alignat}{4}
		\forall \mu \in\eimtr{{a^{\nug}}}, &\qquad& (\eimop{{a^{\nug}}}(\mu))_{\is{{a^{\nug}}}\js{{a^{\nug}}}}&= (\Agamnm(\mu))_{\is{{a^{\nug}}}\js{{a^{\nug}}}},&\qquad &  \forall s \in \indset{1}{\eimcar{{a^{\nug}}}},\label{equa:4:eimcoefA}\\
		\forall (\mu,k) \in 	\eimtr{{b^{\nug}}}, &\qquad & (\eimop{{b^{\nug}}}(\mu,k))_{\is{{b^{\nug}}}\js{{b^{\nug}}}}&=(\Bgamnm(\mu,\ubmka{k}(\mu)))_{\is{{b^{\nug}}}\js{{b^{\nug}}}},&\qquad&  \forall s \in \indset{1}{\eimcar{{b^{\nug}}}}, \label{equa:4:eimcoefB}\\
		\forall (\mu,k) \in \eimtr{{\theta^{\nug}}}, &\qquad & (\eimop{{\theta^{\nug}}}(\mu,k))_{\js{{\theta^{\nug}}}}&=(\Rhonm(\mu,\ubmka{k}(\mu)))_{\js{{\theta^{\nug}}}},&\qquad & \forall s \in \indset{1}{\eimcar{{\theta^{\nug}}}},\label{equa:4:eimcoefTheta} \\
		\forall \mu \in\eimtr{f}, &\qquad& (\eimop{f}(\mu))_{\js{f}}&= (\F(\mu))_{\js{f}},&\qquad &  \forall s \in \indset{1}{\eimcar{f}},\label{equa:4:eimcoefF}
	\end{alignat}
\end{subequations}
and $\eimop{{a^{\nug}}}(\mu)$, $\eimop{{b^{\nug}}}(\mu,k)$, {$\eimop{{\theta^{\nug}}}(\mu,k)$} and {$\eimop{f}(\mu)$} defined in~\eqref{equa:4:EIMapp}.
Notice that we use the HF solution $\ubmka{k}(\mu)$ instead of the RB solution $\ubmredka{k}(\mu)$ to compute the functions $\Alphas{s}{{b^{\nug}}}$ and {$\Alphas{s}{{\theta^{\nug}}}$.}

In the online phase, for every new value of the parameter pair $(\mu, k) \in \mathcal{D} \times \mathbb{N}$, the functions  $\Alphas{s}{{a^{\nug}}}$, $\Alphas{s}{{b^{\nug}}}$, {$\Alphas{s}{{\theta^{\nug}}}$} and {$\Alphas{s}{f}$} are approximated by functions $\AlphaNNs{s}{{a^{\nug}}}$, $\AlphaNNs{s}{{b^{\nug}}}$, {$\AlphaNNs{s}{{\theta^{\nug}}}$}  and {$\AlphaNNs{s}{f}$} which solve the following linear systems:
\begin{subequations}
	\label{equa:4:eimcoefRB}
	\begin{alignat}{5}
		\eimmatrUp{{a^{\nug}}} \AlphaNN{{a^{\nug}}}(\mu) &= \eimsec{{a^{\nug}}}(\mu), &\quad  \eimsec{{a^{\nug}}}(\mu)&:=\big((\Agamnm(\mu))_{\is{{a^{\nug}}}\js{{a^{\nug}}}}\big)_{s \in \indset{1}{\eimcar{{a^{\nug}}}}} \in \Rd{\eimcar{{a^{\nug}}}},\label{equa:4:eimcoefRBA}\\
		\eimmatrUp{{b^{\nug}}} \AlphaNN{{b^{\nug}}}(\mu,k) &= \eimsec{{b^{\nug}}}(\mu,k), &\quad  \eimsec{{b^{\nug}}}(\mu,k)&:=\big((\Bgamnm(\mu,\ubmredka{k}(\mu)))_{\is{{b^{\nug}}}\js{{b^{\nug}}}}\big)_{s \in \indset{1}{\eimcar{{b^{\nug}}}}}\in \Rd{\eimcar{{b^{\nug}}}}, \label{equa:4:eimcoefRBB}\\
		\eimmatrUp{{\theta^{\nug}}} \AlphaNN{{\theta^{\nug}}}(\mu,k) &= \eimsec{{\theta^{\nug}}}(\mu,k), &\quad  \eimsec{{\theta^{\nug}}}(\mu,k)&:=\big((\Rhonm(\mu,\ubmredka{k}(\mu)))_{\js{{\theta^{\nug}}}}\big)_{s \in \indset{1}{\eimcar{{\theta^{\nug}}}}} \in \Rd{\eimcar{{\theta^{\nug}}}},\label{equa:4:eimcoefRBTheta}\\
		\eimmatrUp{f} \AlphaNN{f}(\mu) &= \eimsec{f}(\mu), &\quad  \eimsec{f}(\mu)&:=\big((\F(\mu))_{\js{f}}\big)_{s \in \indset{1}{\eimcar{f}}} \in \Rd{\eimcar{f}},\label{equa:4:eimcoefRBF}
	\end{alignat}
\end{subequations}
where the vector-valued functions $\AlphaNN{{a^{\nug}}}$, $\AlphaNN{{b^{\nug}}}$, $\AlphaNN{{\theta^{\nug}}}$ and $\AlphaNN{f}$ are such that 
\begin{subequations}
	\begin{alignat}{2}
		\AlphaNN{{a^{\nug}}}(\mu)&:=\Vect{\AlphaNNs{s}{{a^{\nug}}}}{\mu}{s}{1}{\eimcar{{a^{\nug}}}} \in \Rd{\eimcar{{a^{\nug}}}},\\ \AlphaNN{{b^{\nug}}}(\mu,k)&:=\Vect{\AlphaNNs{s}{{b^{\nug}}}}{\mu}{s}{1}{\eimcar{{b^{\nug}}}} \in \Rd{\eimcar{{b^{\nug}}}}, \\ \AlphaNN{{\theta^{\nug}}}(\mu,k)&:=\Vect{\AlphaNNs{s}{{\theta^{\nug}}}}{\mu}{s}{1}{\eimcar{{\theta^{\nug}}}} \in \Rd{\eimcar{{\theta^{\nug}}}},\\
		\AlphaNN{f}(\mu)&:=\Vect{\AlphaNNs{s}{f}}{\mu}{s}{1}{\eimcar{f}} \in \Rd{\eimcar{f}}.
	\end{alignat}
\end{subequations}
Notice that here, in the online phase, we use the RB solutions $\ubmredka{k}(\mu)$.
The parameter-independent interpolation matrices $\eimmatrUp{{a^{\nug}}} \in \Rd{\eimcar{{a^{\nug}}} \times \eimcar{{a^{\nug}}}}$, $\eimmatrUp{{b^{\nug}}} \in \Rd{\eimcar{{b^{\nug}}} \times \eimcar{{b^{\nug}}}}$, $\eimmatrUp{{\theta^{\nug}}} \in \Rd{\eimcar{{\theta^{\nug}}} \times \eimcar{{\theta^{\nug}}}}$ and $\eimmatrUp{f} \in \Rd{\eimcar{f} \times \eimcar{f}}$ are such that
\begin{subequations}
	\label{equa:4:eimintmat}
	\begin{alignat}{4}
		(\eimmatrUp{{a^{\nug}}})_{is}&:= (\Agamns)_{\im{{a^{\nug}}}\jm{{a^{\nug}}}},&\qquad &  \forall i,s \in \indset{1}{\eimcar{{a^{\nug}}}},	\label{equa:4:eimintmatA}\\
		(\eimmatrUp{{b^{\nug}}})_{is}&:=(\Bgamns)_{\im{{b^{\nug}}}\jm{{b^{\nug}}}},&\qquad&  \forall i,s \in \indset{1}{\eimcar{{b^{\nug}}}},	\label{equa:4:eimintmatB}\\
		(\eimmatrUp{{\theta^{\nug}}})_{is}&:=(\Rhons)_{\jm{{\theta^{\nug}}}},&\qquad & \forall i,s \in \indset{1}{\eimcar{{\theta^{\nug}}}},	\label{equa:4:eimintmatTheta}\\
		(\eimmatrUp{f})_{is}&:=(\Fs)_{\jm{f}},&\qquad & \forall i,s \in \indset{1}{\eimcar{f}}.	\label{equa:4:eimintmatF}
	\end{alignat}
\end{subequations}
By construction, these matrices are lower-triangular with unit diagonal. Consequently, these matrices are invertible and their inverse can be easily computed once and for all during the offline phase. Combining~\eqref{equa:4:RBmatr} with~\eqref{equa:4:EIMapp}, we obtain the following approximate decompositions:
\begin{subequations}
	\label{equa:4:EIMDECOMPRB}
	\begin{alignat}{4}
		\AgamnM(\mu)&\approx{}&\eimopN{{a^{\nug}}}(\mu):=& \redsum{s}{1}{{a^{\nug}}}\AlphaNNs{s}{{a^{\nug}}}(\mu)\AgamnMs,& \quad \AgamnMs&:= \Z^\top \Agamns \Z \in \mathbb{R}^{N \times N},\label{equa:4:EIMDECOMPRBA}\\
		\BgamnM(\mu,k)&\approx{} &\eimopN{{b^{\nug}}}(\mu,k):=&\redsum{s}{1}{{b^{\nug}}}\AlphaNNs{s}{{b^{\nug}}}(\mu,k)\BgamnMs,& \quad \BgamnMs&:= \Z^\top \Bgamns \Z \in \mathbb{R}^{N \times N},\label{equa:4:EIMDECOMPRBB}\\
		\RhonM(\mu,k)&\approx{}&\eimopN{{\theta^{\nug}}}(\mu,k):=&\redsum{s}{1}{{\theta^{\nug}}}\AlphaNNs{s}{{\theta^{\nug}}}(\mu,k)\RhonMs,& \quad \RhonMs&:= \Z^\top \Rhons \in \mathbb{R}^{N},\label{equa:4:EIMDECOMPRBTheta}\\
		\FM(\mu)&\approx&\eimopN{f}(\mu):=& \redsum{s}{1}{f}\AlphaNNs{s}{f}(\mu)\FMs,& \quad \FMs&:= \Z^\top \Fs \in \mathbb{R}^{N},\label{equa:4:EIMDECOMPRBF}
	\end{alignat}
\end{subequations}
which lead to an efficient offline/online decomposition since the parameter-independent arrays $\AgamnMs$, $\BgamnMs$, {$\RhonMs$} and {$\FMs$} are small-dimensional and can be computed once and for all during the offline phase. Finally, using the approximations from~\eqref{equa:4:EIMDECOMPRB} in~\eqref{equa:4:itNWOF2refred} (for simplicity, we keep the same notation for the unknown $\Deltaubmredka{k}(\mu)$), we consider the following sequence of problems:
For all $(\mu,k) \in \mathcal{D} \times \mathbb{N}$, find $\Deltaubmredka{k}(\mu) \in \mathbb{R}^{\dimPred}$ such that
\begin{equation}
	\eimopN{{a^{\nug}}}(\mu) \Deltaubmredka{k}(\mu) + \eimopN{{b^{\nug}}}(\mu,k) \Deltaubmredka{k}(\mu) = - \eimopN{{r^{\nug}}}(\mu,k),
	\label{equa:4:itNWOF2refredEIM}
\end{equation}
{where $\eimopN{{r^{\nug}}}(\mu,k)$ is given by
\begin{equation}
	\eimopN{{r^{\nug}}}(\mu,k):= \eimopN{{a^{\nug}}}(\mu)\Ubmredka{k}(\mu) + \eimopN{{\theta^{\nug}}}(\mu,k)-\eimopN{f}(\mu).
\end{equation}}%
The convergence criterion for~\eqref{equa:4:itNWOF2refredEIM} is still~\eqref{equa:4:cvcrit} using the reconstructed functions. We denote the converged solution to the sequence of problems~\eqref{equa:4:itNWOF2refredEIM} as $\Ubmredka{\rm{cv}}(\mu) \in \Rd{N}$ and the associated reconstructed solution as $\ubmredka{\rm{cv}}(\mu) \in \fespacebm{V}(\mu)$.

To summarize, our RB procedure is organized as follows:
\begin{itemize}
	\item[---] Offline stage
	\begin{enumerate}
		\item Select a training subset $\discrd{D}{train}:=\set{\mu}{p}{1} \subset \mathcal{D}$.
		\item Compute the HF snapshots $\setbis{\Ubmka{\rm{cv}}(\mup{p})}{p}{1} \subset \Rd{\dimP}$ by solving~\eqref{equa:4:itNWOF2Alg} until convergence on $k$.
		\item Compute the reduced space $\rfespacebm{V}{N}$ by using POD on the snapshots.
		\item Compute the high-dimensional arrays $\Agamns$, $\Bgamns$, {$\Rhons$} and {$\Fs$} by using the EIM.
		\item Invert the interpolation matrices $\eimmatrUp{{a^{\nug}}}$, $\eimmatrUp{{b^{\nug}}}$, {$\eimmatrUp{{\theta^{\nug}}}$} and {$\eimmatrUp{f}$.}
		\item Compute the small-dimensional arrays $\AgamnMs$, $\BgamnMs$, {$\RhonMs$} and {$\FMs$} by using~\eqref{equa:4:EIMDECOMPRB}.
	\end{enumerate}
	\item[---] Online stage:  For any $\mu \in \mathcal{D}\setminus \discrd{D}{train}$,
	\begin{enumerate}
		\item Evaluate $\eimsec{{a^{\nug}}}(\mu)$ {and $\eimsec{f}(\mu)$} and compute $\AlphaNN{{a^{\nug}}}(\mu)$ {and $\AlphaNN{f}(\mu)$} by using~\eqref{equa:4:eimcoefRBA} {and \eqref{equa:4:eimcoefRBF}.}
		\item Evaluate $\eimopN{{a^{\nug}}}(\mu)$ and {$\eimopN{f}(\mu)$} by using~\eqref{equa:4:EIMDECOMPRBA} and {\eqref{equa:4:EIMDECOMPRBF}.}
		\item Loop on $k\geq 0$,
		\begin{enumerate}
			\item \label{item:4:RBalgo2a}Evaluate $\eimsec{{b^{\nug}}}(\mu,k)$ and {$\eimsec{{\theta^{\nug}}}(\mu,k)$} and  compute $\AlphaNN{{b^{\nug}}}(\mu,k)$ and {$\AlphaNN{{\theta^{\nug}}}(\mu,k)$} by using~\eqref{equa:4:eimcoefRBB} and {\eqref{equa:4:eimcoefRBTheta}}.
			\item Evaluate $\eimopN{{b^{\nug}}}(\mu,k)$ and {$\eimopN{{\theta^{\nug}}}(\mu,k)$} by using~\eqref{equa:4:EIMDECOMPRBB} and {\eqref{equa:4:EIMDECOMPRBTheta}}. 
			\item Solve~\eqref{equa:4:itNWOF2refredEIM}. 
			\item Check convergence; if not, set $k \leftarrow k + 1$ and go back to~(\ref{item:4:RBalgo2a}).
		\end{enumerate}
	\end{enumerate}
\end{itemize}

For the contact problem with friction, the RB formulation is obtained in exactly the same way. We simply replace the forms $\Agamn(\mu;\cdot,\cdot)$, $\Bgamn(\mu;\cdot;\cdot,\cdot)$ and $\Resn(\mu;\cdot;\cdot)$ by the forms $\Agam(\mu;\cdot,\cdot)$, $\Bgam(\mu;\cdot;\cdot,\cdot)$ and $\Res(\mu;\cdot;\cdot)$. Notice that for the EIM approximation, additional affine parametric decompositions are performed on the large-dimensional arrays $\Bgamtm(\mu,k)$ and $\Rhotm(\mu,k)$.

\section{Numerical results}\label{section:4:numres}
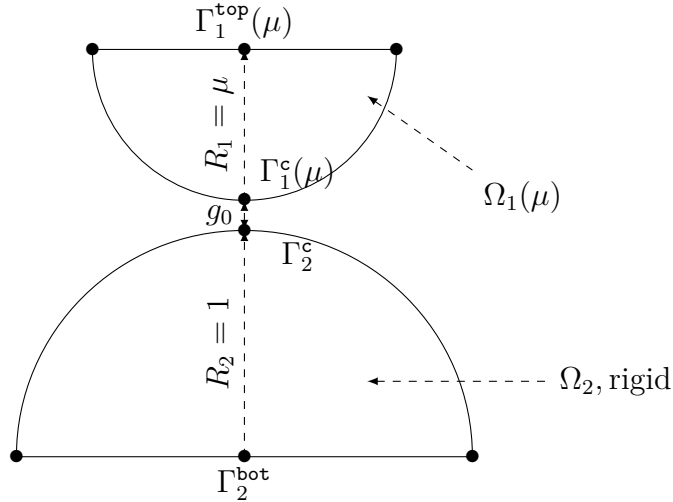
\begin{figure}[!h]
	\centering
	\begin{tikzpicture}[scale = 2]
		%\draw (0,0.2) -- (0,1.2)node[midway,above,sloped] {$\Gamma_{1}^{\texttt{left}}$};
		\draw (-1,1.2) -- (1,1.2)node[midway,above,sloped] {$\Gamma_{1}^{\texttt{top}}(\mu)$};
		\draw (1,1.2) arc (0:-180:1)node[midway,above,sloped,xshift=0.7cm] {$\Gamma_{1}^{\texttt{c}}(\mu)$};
		%\draw  (0,-1.5) -- (0,0) node[midway,above,sloped] {$\Gamma_{2}^{\texttt{left}}$};
		
		\draw (-1.5,-1.5) -- (1.5,-1.5)node[midway,below,sloped] {$\Gamma_{2}^{\texttt{bot}}$};
		\draw (1.5,-1.5) arc (0:180:1.5)node[midway,below,sloped,xshift=0.7cm] {$\Gamma_{2}^{\texttt{c}}$};
		\draw (1,1.2) node{$\bullet$} ;
		\draw (-1,1.2) node{$\bullet$} ;
		\draw (0,1.2) node{$\bullet$} ;
		\draw (0,0.2) node{$\bullet$} ;
		\draw (0,0) node{$\bullet$} ;
		\draw (0,-1.5) node{$\bullet$} ;
		\draw (1.5,-1.5) node{$\bullet$} ;
		\draw (-1.5,-1.5) node{$\bullet$} ;
		\draw [dashed,>=latex,<-] (0.807,0.893) -- (1.5,0.4);
		\draw (1.5,0.4) node[below right] {$\Omega_{1}(\mu)$};
		\draw[dashed,>=latex,<-]  (0.807,-1.) -- (2.,-1.);
		\draw (2.,-1.) node[right] {$\Omega_{2}, \rm{rigid}$};;
		\draw [dashed,>=latex,->] (0,0.2) -- (0,1.2)node[midway,above,sloped] {$R_{1} = \mu$};
		\draw [dashed,>=latex,->] (0,-1.5) -- (0,0)node[midway,above,sloped] {$R_{2} = 1$};
		\draw [dashed,>=latex,<->] (0,0) -- (0,0.2)node[midway,left] {$g_{0}$};
	\end{tikzpicture}
	\caption{Frictionless Hertz test case: Configuration.}
\label{fig:4:HertzConf}
\end{figure}

We consider the Hertz contact problem between the two half-disks as represented in Figure~\ref{fig:4:HertzConf}. The upper half-disk occupies the deformable domain $\Omega_{1}(\mu) \subset \mathbb{R}^2$ of parametric radius
\begin{align}
	R_{1}(\mu):= \mu \quad \mbox{ with  } \mu\in \mathcal{D}:= \big[0.7, 1.3 \big](\SI{}{\metre}),
\end{align} 
and the lower half-disk the rigid domain $\Omega_{2}\subset \mathbb{R}^2$ of fixed radius $R_{2}:=1\SI{}{\metre}$. The initial gap between the two half-disks is equal to $g_{0}>0$. We impose a displacement of $-d$  on $\Gamma_{1}^{\texttt{top}}(\mu)$ of $\Omega_{1}(\mu)$ with  $d \geq g_{0}$. The initial gap $g_0$ and the imposed displacement $d$ are, respectively, set to $g_{0}:= 0.001\SI{}{\metre}$ and $d:= 0.09\SI{}{\metre}$. This latter value, which is less than $10\%$ of the maximum value of $R_1(\mu)$ allows us to remain within the validity of the small deformation assumption. Notice that since $\Omega_{2}$ is rigid and fixed, we only mesh the domain $\Omega_{1}(\mu)$ and set $\Omega(\mu):= \Omega_{1}(\mu)$ to build the HF space. The material parameters are $E:=15\SI{}{\pascal}$ for the Young modulus and $\nu:=0.35$ for the  Poisson coefficient. 
A first training set is typically chosen as $\discrd{D}{train}:=\big\{0.7 + 0.0075i, 0\leq i \leq 60 \big\}(\SI{}{\metre})$ (altogether $P=61$ points), and  the validation set  $\discrd{D}{valid}$ is generated by choosing $30$ elements in $\mathcal{D}$ randomly with a uniform distribution. A second richer training set (with $P = 201$ points) will also be considered.

We consider the reference domain $\omegaahat_{1}:= \Omega_1(1)$ and introduce the geometric mapping $h_{1}(\mu): \omegaahat_{1} \rightarrow \Omega_{1}(\mu)$ defined as $h_{1}(\mu)(\xbm):=\mu\xbm$, for all $\xbm \in  \omegaahat_{1}$, with the origin located at the center of $\omegaahat_{1}$. We use a mesh composed of ${1956}$ nodes with ${633}$ nodes on the potential contact manifold $\gammaachat_{1}$ which is the part of the half circle $\gammaachat_{1}$ of angle $\theta \in [-\frac{5\pi}{8},-\frac{3\pi}{8}]$ with respect to the horizontal axis. For all $\mu \in \mathcal{D}$, we equip the space $\fespacebm{V}(\mu)$ with the norm $\norm{\cdot}{\fespacebm{V}(\mu)}$ defined as follows: 
\begin{align}
	\norm{\vbm}{\fespacebm{V}(\mu)}:= \Big( \norm{\vbm}{L^2(\omegaa(\mu))}^2 + \widehat{\ell}^2 \norm{\nabla \vbm}{L^2(\omegaa(\mu))}^2\Big)^{\frac{1}{2}},\quad \forall \vbm \in \fespacebm{V}(\mu),
\end{align}
where the characteristic length $\widehat{\ell}:=1$ is the radius of $\omegaahat_{1}$ and is introduced for dimensional consistency. The HF and RB computations use the python library of the finite element software \getfem~\cite{Renard-2020}.

\subsection{Frictionless case}\label{subsection:4:frictionless}
We first consider the frictionless Hertz contact problem.
\subsubsection{Results using Nitsche's method}\label{subsubsection:4:N}
For the discretization, we use $\mathbb{P}_{2}$ Lagrange finite elements leading to ${\dimP:=14918}$ degrees of freedom. We choose $\gamma:=\frac{\gamma_{0}}{h}$ with {$\gamma_{0}:= 50\mulame$}, where ${h:=2.5\SI{}{\milli\metre}}$ refers to the mesh size and $\mulame:= \frac{E}{2(1 + \nu)}\SI{}{\pascal}$ refers to the second Lamé parameter.

\begin{figure}[!h]
	\centering
	\begin{minipage}[l]{.45\linewidth}
		\includegraphics[width=6.5cm]{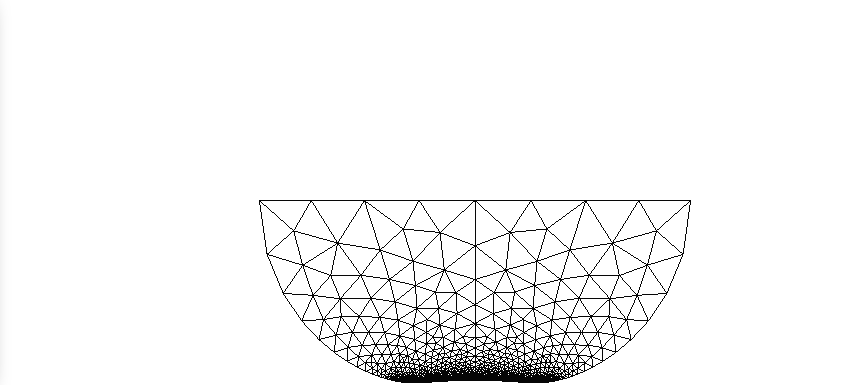}
	\end{minipage} 
	\begin{minipage}[l]{.45\linewidth}
		\includegraphics[width=6.5cm]{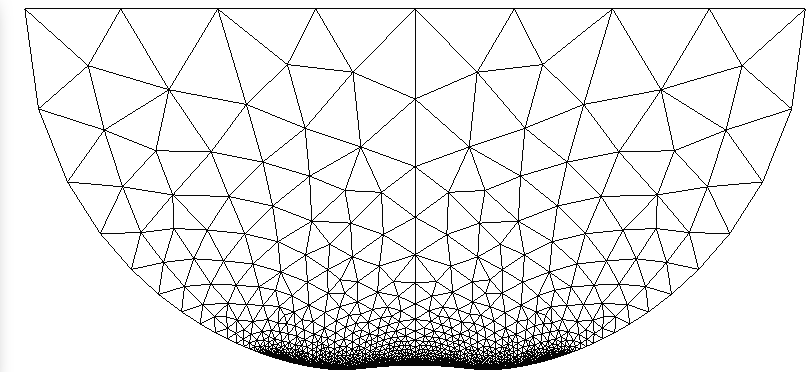}
	\end{minipage}
	\caption {Frictionless Hertz test case, Nitsche's method (HF): $\ubm(\mu)$. Left: $\mu= 0.7 \SI{}{\metre}$; Right: $\mu= 1.3 \SI{}{\metre}$.}
	\label{fig:4:HertzLN:defconf}
\end{figure}

Figure~\ref{fig:4:HertzLN:defconf} displays the deformed configurations resulting from the HF displacement fields $\ubm(\mu)$ for $\mu= 0.7 \SI{}{\metre}$ (left panel) and for $\mu=  1.3\SI{}{\metre}$ (right panel). {We can see that we use a symmetric mesh. This is important because it guarantees the symmetry of the HF snapshots. Indeed, if the snapshots are not  symmetric, the resulting POD modes will not be either. Consequently, the reduced model looses this symmetry property, leading to reduced solutions of poorer quality. Moreover, we have discretized a complete half-disk instead of a quarter-disk to avoid some difficulties when enforcing the symmetry condition at the lowest point of $\omegaa_{1}(\mu)$ in the case of the frictional contact problems (see Section \ref{subsection:4:FC}).}

\begin{figure}[!h]
	\centering
	\begin{subfigure}{0.32\textwidth}
		\includegraphics[width=0.8\textwidth]{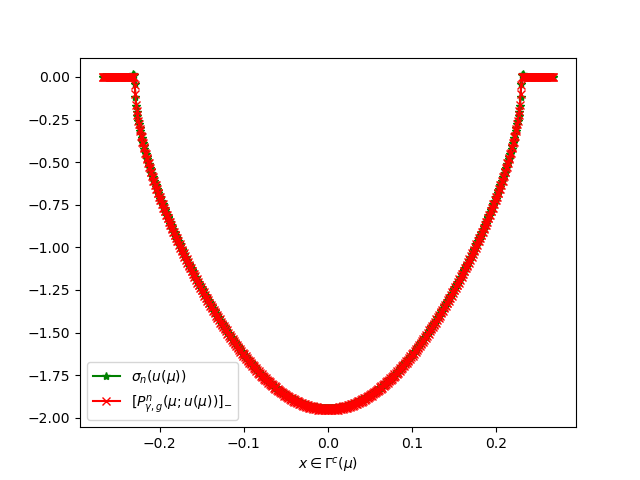}
		%\caption{$\sigman(\ubm(\mu))$.}
		\label{fig:4:HertzLN:sigma07}
	\end{subfigure}
	\hfill
	\begin{subfigure}{0.32\textwidth}
		\includegraphics[width=0.8\textwidth]{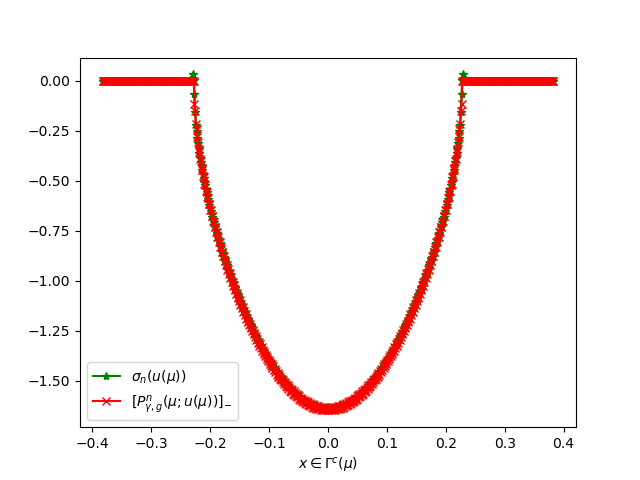}
		%\caption{$\un(\mu) - g(\mu)$.}
		\label{fig:4:HertzLN:sigma10}
	\end{subfigure}
		\hfill
	\begin{subfigure}{0.32\textwidth}
		\includegraphics[width=0.8\textwidth]{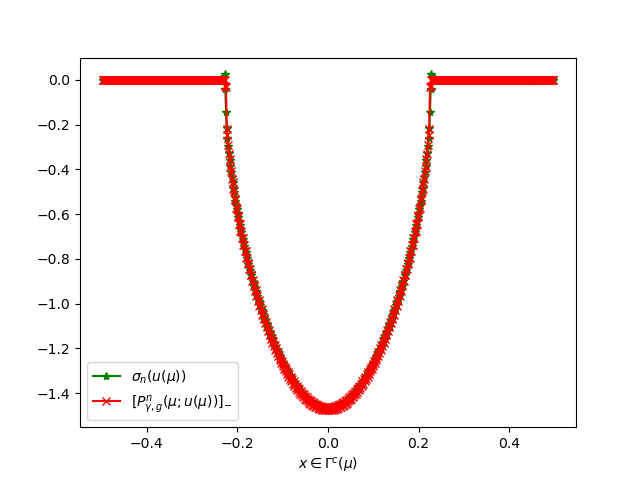}
		%\caption{$\un(\mu) - g(\mu)$.}
		\label{fig:4:HertzLN:sigma13}
	\end{subfigure}

	\begin{subfigure}{0.32\textwidth}
		\includegraphics[width=0.8\textwidth]{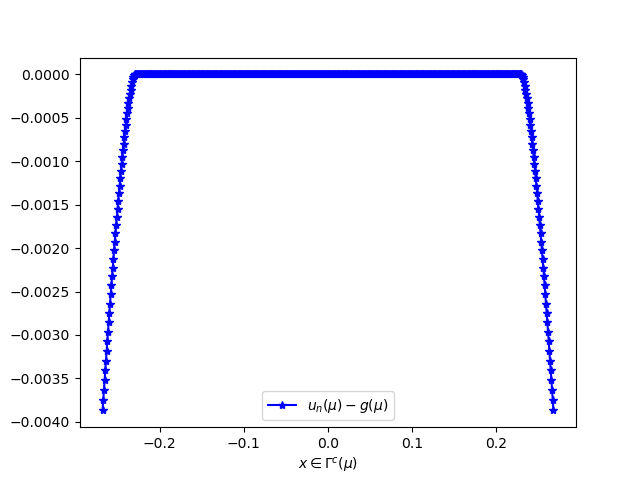}
		%\caption{$\sigman(\ubm(\mu))$.}
		\label{fig:4:HertzLN:gap07}
	\end{subfigure}
	\hfill
	\begin{subfigure}{0.32\textwidth}
		\includegraphics[width=0.8\textwidth]{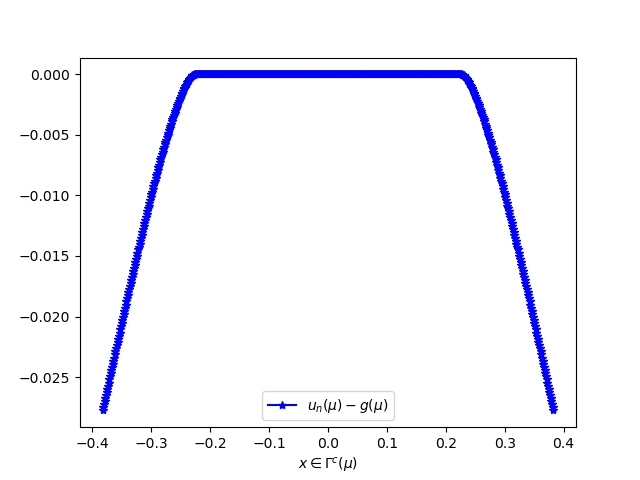}
		%\caption{$\un(\mu) - g(\mu)$.}
		\label{fig:4:HertzLN:gap10}
	\end{subfigure}
		\hfill
	\begin{subfigure}{0.32\textwidth}
		\includegraphics[width=0.8\textwidth]{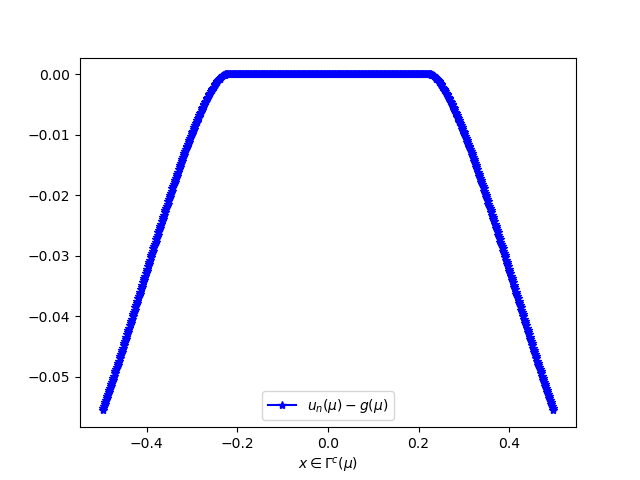}
		%\caption{$\un(\mu) - g(\mu)$.}
		\label{fig:4:HertzLN:gap13}
	\end{subfigure}
	
	\caption{Frictionless Hertz test case, Nitsche's method (HF).  Row~$1$: $\sigman(\ubm(\mu))$ and $\projconvb{\Pgamn(\mu;\ubm(\mu))}{-}$. Row~$2$: $\un(\mu) - g(\mu)$. Column~$1$: $\mu=0.7\SI{}{\metre}$. Column~$2$: $\mu=1.0\SI{}{\metre}$. Column~$3$: $\mu=1.3\SI{}{\metre}$.}
	\label{fig:4:HertzLN:CU-AL}
\end{figure}

{Figure~\ref{fig:4:HertzLN:CU-AL} displays the superposition of the normal stress $\sigman(\ubm(\mu))$ and its Alart--Curnier reformulation $\projconv{\Pgamn(\mu;\ubm(\mu))}{-}$  (Row~$1$) and the gap on the deformed configuration $\un(\mu) - g(\mu)$ (Row~$2$) as a function of the abscissa along $\gammaac(\mu)$ for $\mu = 0.7\SI{}{\metre}$ (Column~$1$), $\mu = 1\SI{}{\metre}$ (Column~$2$), and $\mu = 1.3\SI{}{\metre}$ (Column~$3$).}
We see that the normal stress is of good quality (with almost no spurious oscillations) and matches very well with its counterpart $\projconv{\Pgamn(\mu;\ubm(\mu))}{-}$ resulting from the Alart--Curnier reformulation. We also see that the negativity condition on the gap $\un(\mu) - g(\mu)$ is satisfied on the whole potential contact manifold $\gammaac(\mu)$.
{To have a better look at Signorini's contact conditions, we display in Table~\ref{table:4:HertzLN:signo} the relative error on the Alart--Curnier reformulation of Signorini's contact conditions defined as follows:
\begin{equation}
e_{\rm{AC}}^{\nug}(\mu):= \frac{\norm{\sigman(\ubm(\mu)) - \projconvb{\Pgamn(\mu;\ubm(\mu))}{-}}{\ell^{2}(\gammaac(\mu))}}{\norm{\sigman(\ubm(\mu))}{\ell^{2}(\gammaac(\mu))}},
\end{equation}
where the discrete $\ell^{2}(\gammaac(\mu))$-norms are sampled at the mesh nodes located on $\gammaac(\mu)$.
For the three values of the parameter $\mu \in \{0.7,1,1.3\}(\SI{}{\metre})$, we consider three values of the mesh size $h$, namely $h=5\SI{}{\milli\metre}$ (coarse), $h=2.5\SI{}{\milli\metre}$ (medium) and $h=1.25\SI{}{\milli\metre}$ (fine). We notice that the relative error $e_{\rm{AC}}^{\nug}(\mu)$ is smaller than $1.5\%$ for the three parameters values and for the three mesh sizes. Moreover, we see that the error decreases when $h$ decreases, thereby indicating the convergence of the approximation. More precisely, we observe a convergence of order $1$ for the approximation of $\sigman(\ubm(\mu))$ (with a slight order reduction for the larger value of $\mu$). Thus, we can say that the Signorini contact conditions are globally satisfied with a good accuracy although they are not strictly enforced.}
\begin{table}[!h]
	\centering
	\renewcommand{\arraystretch}{1.5}
	\begin{tabular}{|c|*{9}{>{\centering\arraybackslash}m{1.1cm}|} }
		\hline
		$\mu(\SI{}{\metre})$ & \multicolumn{3}{|c|}{$0.7$}&  \multicolumn{3}{|c|}{$1$}&  \multicolumn{3}{|c|}{$1.3$} \\
		\hline
		$h(\SI{}{\milli\metre})$&$5$ & $2.5$& $1.25$ &$5$ & $2.5$& $1.25$ &$5$ & $2.5$& $1.25$ \\
		\hline
		$e_{\rm{AC}}^{\nug}(\%)$ &$1$& $0.52$& $0.28$& $1.45$ &$0.72$& $0.33$& $1.17$& $1.1$&$0.43$ \\
		\hline
	\end{tabular}
	\caption{{Frictionless Hertz test case, Nitsche's method (HF): Relative error $e_{\rm{AC}}^{\nug}(\mu)$ for $\mu \in \{0.7,1.0,1.3\}(\SI{}{\metre})$ and the mesh sizes $h \in \{5,2.5,1.25\}(\SI{}{\milli\metre})$.}}
	\label{table:4:HertzLN:signo}
\end{table}

\begin{figure}[!h]
	\centering
	\includegraphics[width=8cm]{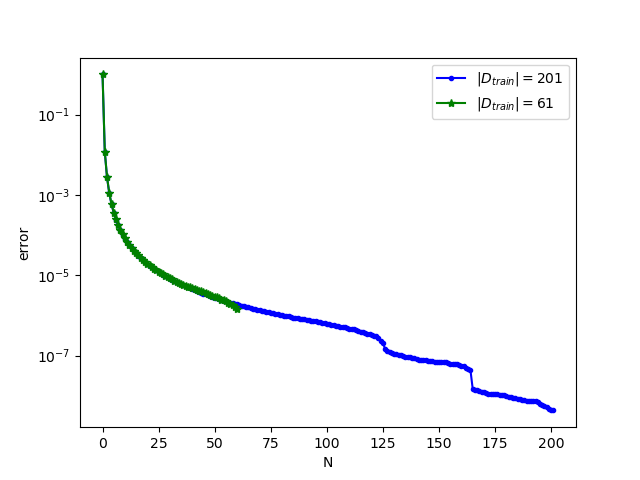}
	\caption {Frictionless Hertz test case, Nitsche's method: {Relative} POD projection error {$e_{\texttt{POD}}(N)$ as function of $\dimPred$, for $|\discrd{D}{train}|=61$ and ${|\discrd{D}{train}|=201}$.}}
	\label{fig:4:HertzLN:POD_error}
\end{figure}

Let us consider the relative POD projection error defined as follows:
\begin{align}
	e_{\texttt{POD}}(N):=\frac{ \Big(\sum\limits_{p \in\indset{1}{P} }  \norm{ \big(\mathbb{I}_{\Rd{\dimP}}- \proj{\rfespacebm{V}{\dimPred}}{}\big)( \Ubmka{\rm{cv}}(\mup{p})) }{{\stiff(\mu_{p})}}^2\Big)^{\frac{1}{2}}}{ \Big(\sum\limits_{p \in\indset{1}{P} }  \norm{\Ubmka{\rm{cv}}(\mup{p})}{{\stiff(\mu_{p})} }^2\Big)^{\frac{1}{2}}},
	\label{equa:4:PODerror}
\end{align}
where $\proj{\rfespacebm{V}{\dimPred}}{}$ denotes the orthogonal projection onto $\rfespacebm{V}{\dimPred} \subset \Rd{\dimP}$ and  {$\stiff(\mu_{p})$} the Gram matrix of the inner product associated with {$H^{1}(\omegaa(\mu_{p});\Rd{d})$}. We consider two different training sets, the first with  cardinality $61$ and the second with cardinality {$201$}. Figure~\ref{fig:4:HertzLN:POD_error} shows the relative projection error $e_{\texttt{POD}}(N)$ produced by the POD algorithm as a function of the number of vectors composing the reduced basis {for both training sets}. In all cases, we notice that the projection error decreases sufficiently fast so that indeed the linear spaces generated by the snapshots can be approximated by small-dimensional subspaces. We also observe a fast decrease of the POD error for the first $15$ modes before a slower decrease occurs at error levels between {$10^{-5}$ and $10^{-6}$ (resp. $10^{-5}$ and $10^{-7}$) for the first (resp. second) training set}.

\begin{figure}[!h]
	\centering
	\begin{minipage}[l]{.45\linewidth}
		\includegraphics[width=7cm]{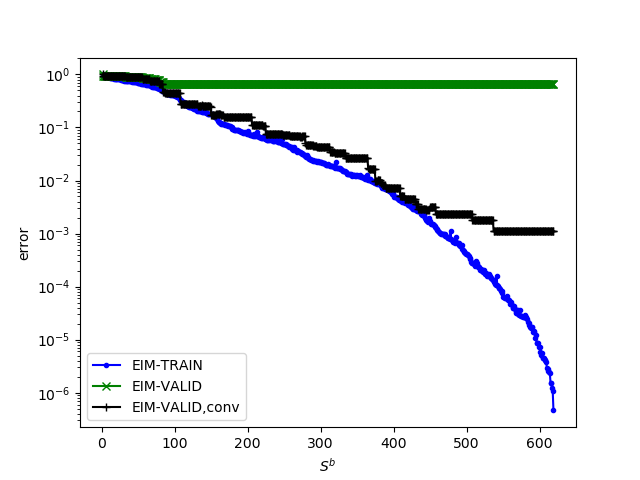}
	\end{minipage} 
	\begin{minipage}[l]{.45\linewidth}
		\includegraphics[width=7cm]{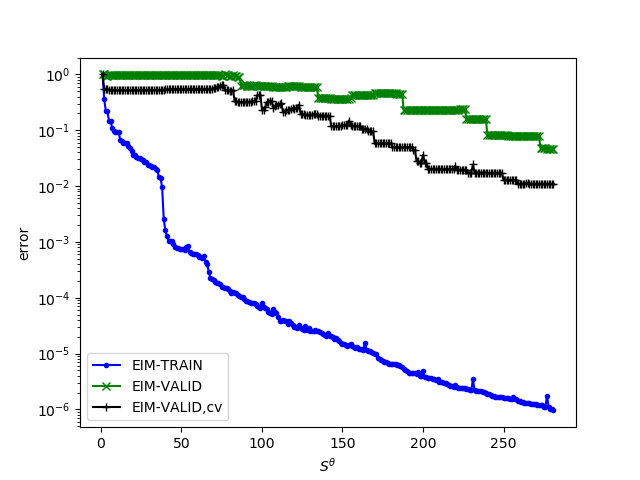}
	\end{minipage}
	\caption {{Frictionless Hertz test case, Nitsche's method: {Relative} EIM approximation errors as a function of the rank $\eimcar{{b^{\nug}}}$ or $\eimcar{{\theta^{\nug}}}$ of the approximation. Left: $e_{\texttt{EIM}}^{{b^{\nug}}}(\eimcar{{b^{\nug}}},\discrd{D}{train})$,  $e_{\texttt{EIM}}^{{b^{\nug}}}(\eimcar{{b^{\nug}}},\discrd{D}{valid})$ and $e_{\texttt{EIM}}^{b^{\nug},\rm{cv}}(\eimcar{{b^{\nug}}})$; Right: $e_{\texttt{EIM}}^{{\theta^{\nug}}}(\eimcar{{\theta^{\nug}}},\discrd{D}{train})$, $e_{\texttt{EIM}}^{{\theta^{\nug}}}(\eimcar{{\theta^{\nug}}},\discrd{D}{valid})$ and $e_{\texttt{EIM}}^{\theta^{\nug},\rm{cv}}(\eimcar{{\theta^{\nug}}})$.}}
	\label{fig:4:HertzLN:eim_error}
\end{figure}

For the EIM approximation, we use the training set $\discrd{D}{train}$ of cardinality 61, and the training sets $\eimtr{{b^{\nug}}}$ and $\eimtr{{\theta^{\nug}}}$ introduced in \eqref{equa:4:eimtrain} are then of cardinality $897$. We fix a tolerance $\tole{EIM}:=10^{-6}$ to bound the errors resulting from~\eqref{equa:4:EIMapp}. With this choice, we obtain {$\eimcar{{b^{\nug}}}=619\ll \dimP\times \dimP$ and $\eimcar{{\theta^{\nug}}}=281 \ll \dimP$}. Notice that in the present test case, we do not need to perform an EIM decomposition on $\Agamnm(\mu)$ since this matrix is already linearly dependent on $\mu$; moreover, {$\F(\mu)$ vanishes} {since we only use an imposed displacement for the load}.
Figure~\ref{fig:4:HertzLN:eim_error} shows the {relative} EIM interpolation errors for the tangent matrix $\Bgamnm(\mu,\ubmka{k}(\mu))$ (left panel) and  the residual vector {$\Rhonm(\mu,\ubmka{k}(\mu))$} (right panel) {as a function of the rank $\eimcar{{b^{\nug}}}$ or $\eimcar{{\theta^{\nug}}}$}, \ie, we plot 
\begin{subequations}
	\begin{alignat}{2}
		e_{\texttt{EIM}}^{{b^{\nug}}}(\eimcar{{b^{\nug}}},\discrd{D}{*})&:= \frac{\max\limits_{\mu \in \discrd{D}{*}}\  \max\limits_{k \in \indset{1}{\kc(\mu)}} \norm{\Bgamnm(\mu,\ubmka{k}(\mu)) - \eimop{{b^{\nug}}}(\mu,k)}{\ell^{\infty}(ij)}}{\max\limits_{\mu \in \discrd{D}{*}}\ \max\limits_{k \in \indset{1}{\kc(\mu)}} \norm{\Bgamnm(\mu,\ubmka{k}(\mu))}{\ell^{\infty}(ij)}}, \label{equa:4:erroreimbn}\\
		e_{\texttt{EIM}}^{{\theta^{\nug}}}(\eimcar{{\theta^{\nug}}},\discrd{D}{*})&:= \frac{\max\limits_{\mu \in \discrd{D}{*}}\ \max\limits_{k \in \indset{1}{\kc(\mu)}} \norm{\Rhonm(\mu,\ubmka{k}(\mu)) - \eimop{{\theta^{\nug}}}(\mu,k)}{\ell^{\infty}(j)}}{\max\limits_{\mu \in \discrd{D}{*}}\ \max\limits_{k \in \indset{1}{\kc(\mu)}} \norm{\Rhonm(\mu,\ubmka{k}(\mu)) }{\ell^{\infty}(j)}},
		\label{equa:4:erroreimthetan}
	\end{alignat}
\end{subequations}
where for a generic matrix (resp.~vector) $M \in \Rd{\dimP\times \dimP}$ (resp.~$V \in \Rd{\dimP}$),
\begin{subequations}
	\begin{alignat}{2}
		\norm{M}{\ell^{\infty}(ij)}&:= \max\limits_{(i,j) \in \indset{1}{\dimP} \times \indset{1}{\dimP}} |M_{ij}|, \\
		\norm{V}{\ell^{\infty}(j)}&:= \max\limits_{j \in \indset{1}{\dimP}} |V_{j}|,
	\end{alignat}
\end{subequations}
{and with $\mathcal{D}_{*}$ either equal to $\discrd{D}{train}$ or to $\discrd{D}{valid}$. Let us first consider the training set $\discrd{D}{train}$. We observe that both errors decrease fast enough to allow accurate approximations. For the tangent matrix $\Bgamnm(\mu,\ubmka{k}(\mu))$, we notice a quasi-uniform decrease of the relative error $e_{\texttt{EIM}}^{{b^{\nug}}}(\eimcar{{b^{\nug}}},\discrd{D}{train})$ with an acceleration at error levels between $10^{-3}$ and $10^{-6}$. For the residual vector $\Rhonm(\mu,\ubmka{k}(\mu))$, we observe a fast decrease of the relative error $e_{\texttt{EIM}}^{{\theta^{\nug}}}(\eimcar{{\theta^{\nug}}},\discrd{D}{train})$ for ranks $\eimcar{{\theta^{\nug}}}$ between $1$ and $40$ yielding errors between $1$ down to $10^{-2}$, and then a significant drop of the error at about $\eimcar{{\theta^{\nug}}}=40$ before a slower decrease occurs at error levels between $10^{-3}$ and $10^{-6}$. Considering the validation set $\discrd{D}{valid}$, we additionally plot the relative EIM approximation errors $e_{\texttt{EIM}}^{b^{\nug},\rm{cv}}(\eimcar{b^{\nug}})$ and $e_{\texttt{EIM}}^{\thetan,\rm{cv}}(\eimcar{\thetan})$ defined  as}
\begin{subequations}
	\begin{alignat}{2}
		e_{\texttt{EIM}}^{b^{\nug},\rm{cv}}(\eimcar{{b^{\nug}}})&:= \frac{\max\limits_{\mu \in \discrd{D}{valid}}  \norm{\Bgamnm(\mu,\ubmka{\rm{cv}}(\mu)) - \eimop{{b^{\nug}}}(\mu,\kc(\mu))}{\ell^{\infty}(ij)}}{\max\limits_{\mu \in \discrd{D}{valid}} \norm{\Bgamnm(\mu,\ubmka{\rm{cv}}(\mu))}{\ell^{\infty}(ij)}}, \label{equa:4:erroreimbnconv}\\
		e_{\texttt{EIM}}^{\theta^{\nug},\rm{cv}}(\eimcar{{\theta^{\nug}}})&:= \frac{\max\limits_{\mu \in \discrd{D}{valid}} \norm{\Rhonm(\mu,\ubmka{\rm{cv}}(\mu)) - \eimop{{\theta^{\nug}}}(\mu,\kc(\mu))}{\ell^{\infty}(j)}}{\max\limits_{\mu \in \discrd{D}{valid}} \norm{\Rhonm(\mu,\ubmka{\rm{cv}}(\mu)) }{\ell^{\infty}(j)}}. \label{equa:4:erroreimthetanconv}
	\end{alignat}
\end{subequations}
{These errors correspond to the {relative} EIM approximation error at convergence of the iterative algorithm, \ie, when $k=\kc(\mu)$. For the tangent matrix $\Bgamnm(\mu,\ubmka{k}(\mu))$, considering first $e_{\texttt{EIM}}^{{b^{\nug}}}(\eimcar{{b^{\nug}}},\discrd{D}{valid})$, we notice a quite modest decrease of the error for ranks $\eimcar{{b^{\nug}}}$ between $1$ and $100$ with errors values between $0.95$ and $0.7$, and then a stagnation of the error at around $0.65$. Considering $e_{\texttt{EIM}}^{b^{\nug},\rm{cv}}(\eimcar{{b^{\nug}}})$, we observe instead a rather uniform decrease of the error, with values quite close to those of $e_{\texttt{EIM}}^{{b^{\nug}}}(\eimcar{{b^{\nug}}},\discrd{D}{train})$, before a stagnation occurs at a value of about $10^{-3}$. For the residual vector $\Rhonm(\mu,\ubmka{k}(\mu))$, considering first $e_{\texttt{EIM}}^{{\theta^{\nug}}}(\eimcar{{\theta^{\nug}}},\discrd{D}{valid})$, we observe a stagnation of the error for ranks $\eimcar{{\theta^{\nug}}}$  between $1$ and $80$ with error values between $1$ and $0.9$, and then  a slower decrease at error levels between $0.8$ and $5\cdot10^{-2}$ with some stagnation phases. Considering $e_{\texttt{EIM}}^{\theta^{\nug},\rm{cv}}(\eimcar{{\theta^{\nug}}})$, we instead observe a stagnation at about $0.5$ for ranks $\eimcar{{\theta^{\nug}}}$ between $2$ and $80$, and then a slower decrease with error values between $0.3$ down to $10^{-2}$. We conclude that the EIM approximation is not very accurate for $\mu \in \discrd{D}{valid}$ and small values of $k$, whereas the accuracy becomes more satisfactory as $k \rightarrow \kc(\mu)$. Therefore, we may expect some difficulties in achieving convergence in the iterative solvers applied to reduced problems, but if convergence is indeed achieved, the accuracy should be reasonable.}

\begin{figure}[!h]
	\centering
	\begin{minipage}[l]{.45\linewidth}
		\includegraphics[width=7cm]{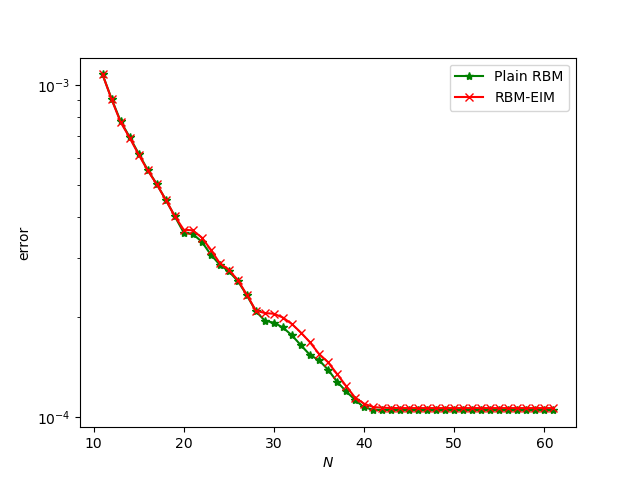}
	\end{minipage} 
	\begin{minipage}[l]{.45\linewidth}
		\includegraphics[width=7cm]{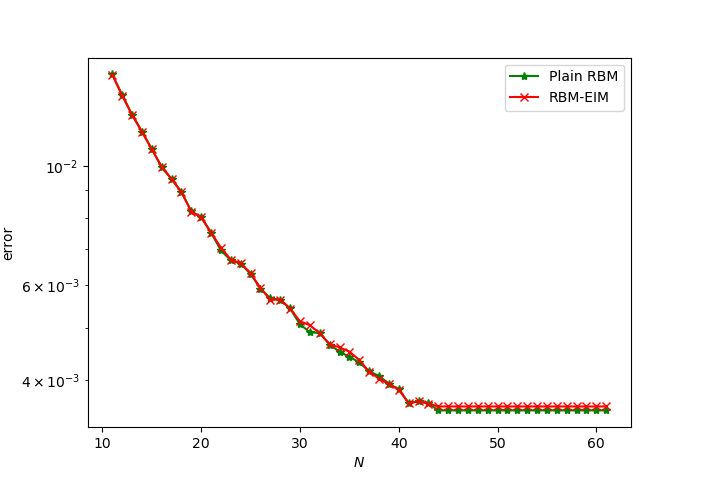}
	\end{minipage}
		\caption {{Frictionless Hertz test case, Nitsche's method: RBM approximation errors  for $|\discrd{D}{train}|=61$. Left: $\e_{\dimPred,\rm{max}}^{\ubm}$; Right: $\e_{\dimPred,\rm{max}}^{\nug\nug}$.}}
	\label{fig:4:HertzLN:RB_error}
\end{figure}

{We denote by $\erb{\ubm}(\mu)$ (resp. $\erb{\nug\nug}(\mu)$) the relative RB approximation error on the displacement (resp. normal stress) defined as}
\begin{subequations}
\begin{alignat}{2}
	\erb{\ubm}(\mu)&:=\frac{\norm{\ubmka{\rm{cv}}(\mu) - \ubmredka{\rm{cv}}(\mu)}{\fespacebm{V}(\mu)}}{\norm{\ubmka{\rm{cv}}(\mu) }{\fespacebm{V}(\mu)}},\\
	\erb{\nug\nug}(\mu)&:=\frac{\norm{\sigman(\ubmka{\rm{cv}}(\mu)) - \sigman(\ubmredka{\rm{cv}}(\mu))}{\ell^2(\gammaac(\mu))}}{\norm{\sigman(\ubmka{\rm{cv}}(\mu)) }{\ell^2(\gammaac(\mu))}},
\end{alignat}
\end{subequations}
and introduce the relative error measures $\e_{\dimPred,\rm{max}}^{u}$ and $\e_{\dimPred,\rm{max}}^{\nug\nug}$ defined as
\begin{equation}
	\e_{\dimPred,\rm{max}}^{\ubm}:= \max\limits_{ \mu\in \discrd{D}{valid}} \erb{\ubm}(\mu),\qquad 	\e_{\dimPred,\rm{max}}^{\nug\nug}:= \max\limits_{ \mu\in \discrd{D}{valid}} \erb{\nug\nug}(\mu).
\end{equation}
Figure~\ref{fig:4:HertzLN:RB_error} displays $\e_{\dimPred,\rm{max}}^{\ubm}$ (left panel) and $\e_{\dimPred,\rm{max}}^{\nug\nug}$ (right panel) as a function of the number of vectors composing the reduced basis. We only consider RB dimensions $N$ larger than $10$. Indeed, for smaller values, the iterative algorithm does not converge for some values of the parameter $\mu$. This can be explained by the poor quality of the EIM approximation of the tangent matrix $\BgamnM(\mu,k)$ for small values of $N$ due to the fact that the RB solution at the first iterations of the iterative algorithm on the reduced model is quite far from the HF solution on which the training of the EIM is performed. We also observe some convergence difficulties for values of the parameter $\mu$ larger than $1.2\SI{}{\metre}$. For this reason, we consider a validation set  $\discrd{D}{valid}$ restricted to the interval $[0.7,1.18](\SI{}{\metre})$. In Figure~\ref{fig:4:HertzLN:RB_error}, we superpose the plain RBM approximation error (without any EIM, thus computationally inefficient) and the RBM-EIM approximation error. We observe that similar errors are obtained for plain RBM and RBM-EIM. This confirms the good quality of the EIM approximations at convergence as claimed above. With both approaches, we notice a stagnation of the relative error at about $10^{-4}$ from about $40$ modes, in agreement with the stagnation observed on the POD projection error.
\begin{comment}
	\begin{figure}[!h]
		\centering
		\includegraphics[width=8cm]{errRBP-CE2.png}
		\caption {Frictionless Hertz test case, Nitsche's method: RBM approximation errors  for $|\discrd{D}{train}|=61$. Left: $\e_{\dimPred,\rm{max}}^{\ubm}$; Right: $\e_{\dimPred,\rm{max}}^{\nug}$.}
		\label{fig:4:HertzLN:RB_error}
	\end{figure}
\end{comment}
\subsubsection{Comparison with the mixed formulation}
\begin{figure}[!h]
	\centering
	\includegraphics[width=8cm]{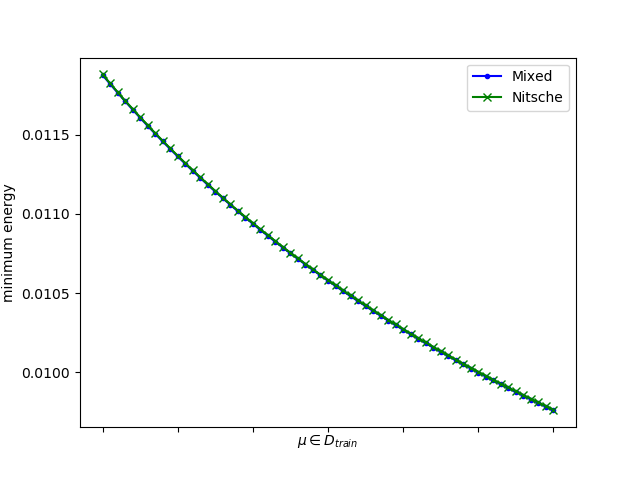}
	\caption {Frictionless Hertz test case: Comparison of the HF energy $\mathcal{J}(\mu;\ubm(\mu))$ between the mixed formulation and Nitsche's method for all $\mu \in \discrd{D}{train}$.}
	\label{fig:4:HertzLM:minenerg}
\end{figure}

For the comparison, we consider the primal-dual formulation employed in~\cite{Niakh-2022} with $\mathbb{P}_{2}$ finite elements for displacement and $\mathbb{P}_{1}$ finite elements for the  Lagrange multiplier. With this choice of the discretization, we can compare on a fair basis the HF displacements obtained with the mixed formulation and with Nitsche's method. We display in Figure~\ref{fig:4:HertzLM:minenerg} the  HF energy $\mathcal{J}(\mu;\ubm(\mu))$ (see~\eqref{equa:4:energy}) for all $\mu \in \discrd{D}{train}$. We notice that we obtain (in the eyeball norm) the same values for the two methods for all $\mu \in \discrd{D}{train}$. Thus, although the constraints are not strictly imposed with Nitsche's method, we obtain a satisfactory accuracy for the quality of the solution in comparison with the mixed formulation. We can also see that the energy decreases with the parameter radius $\mu$ of the half-disk $\omegaa(\mu)$.
The complementary condition is not reported since it is on the order of the machine precision ($10^{-14}\SI{}{\newton/\metre}$) as expected with the primal-dual formulation since the constraints are exactly enforced.

\begin{figure}[!h]
	\centering
	\begin{minipage}[l]{.45\linewidth}
		\includegraphics[width=7cm]{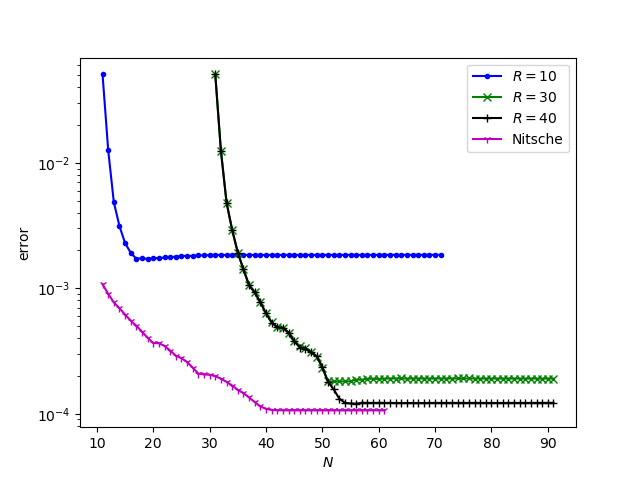}
	\end{minipage} 
	\begin{minipage}[l]{.45\linewidth}
		\includegraphics[width=7cm]{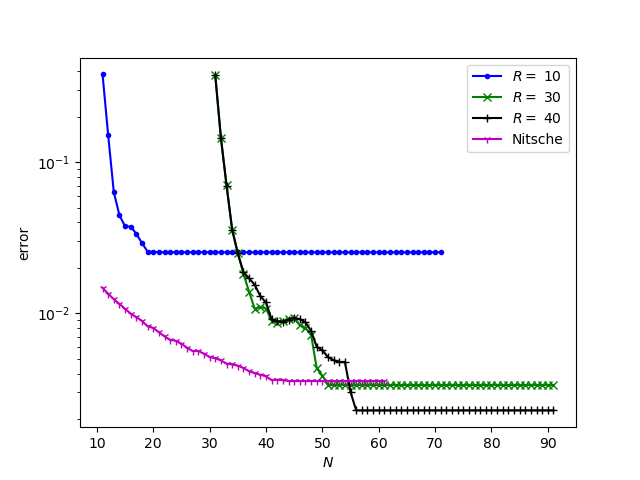}
	\end{minipage}
	\caption {{Frictionless Hertz test case, mixed method: RBM approximation errors. Right: displacement errors $\e_{\dimPred,\dimDred,\rm{max}}^{\ubm}$ and $\e_{\dimPred,\rm{max}}^{\ubm}$; Left: normal stress errors $\e_{\dimPred,\dimDred,\rm{max}}^{\lambda}$ and $\e_{\dimPred,\rm{max}}^{\nug\nug}$.}}
	\label{fig:4:HertzLM:errRBM}
\end{figure}

In contrast to Nitsche's method, in the mixed formulation, it is necessary to stabilize the RBM in order to ensure inf-sup stability for the pair of primal/dual reduced spaces. {For this purpose}, we use the Projected Greedy Algorithm (PGA) algorithm from~\cite{Niakh-2022} with a tolerance $\tole{PGA} := \inf\limits_{\mu \in \discrd{D}{train}} \frac{\betahf(\mu)}{\chf(\mu)}=0.047$ so that the stability condition established in~\cite[Prop~3.1]{Niakh-2022} is fulfilled (the quantities $\betahf(\mu)$ and $\chf(\mu)$ are defined therein). We denote by $\erbpr(\mu)$ (resp.~$\erbdu(\mu)$) the primal (resp.~dual) relative RB approximation error defined as follows:
\begin{align}
	\erbpr(\mu):=\frac{\norm{\ubm(\mu) - \ubmred(\mu)}{\fespacebm{V}(\mu)}}{\norm{\ubm(\mu) }{\fespacebm{V}(\mu)}}, \quad
	\erbdu(\mu):=\frac{\norm{\lambda(\mu) - \lambdared(\mu)}{{\ell^2(\gammaac(\mu))}}}{\norm{\lambda(\mu) }{{\ell^2(\gammaac(\mu))}}},
\end{align}
and introduce the {relative} primal (resp.~dual) error measure $\e_{\dimPred,\dimDred,\rm{max}}^{\ubm}$ (resp.~$\e_{\dimPred,\dimDred,\rm{max}}^{\lambda}$) defined as
\begin{align}
	\e_{\dimPred,\dimDred,\rm{max}}^{\ubm}:= \max\limits_{ \mu\in \discrd{D}{valid}} \erbpr(\mu), \quad
	\e_{\dimPred,\dimDred,\rm{max}}^{\lambda}:= \max\limits_{ \mu\in \discrd{D}{valid}} \erbdu(\mu).
\end{align}
Figure~\ref{fig:4:HertzLM:errRBM} displays the quantity $\e_{\dimPred,\dimDred,\rm{max}}^{\ubm}$ (resp.~ $\e_{\dimPred,\dimDred,\rm{max}}^{\lambda}$) on the left (resp.~right) panel as a function of the dimension of the reduced primal basis (after stabilization) for three values of the dimension of the reduced dual basis $R$, namely $R=10$, $R=30$ and $R=40$. {In addition, we plot the relative error $\e_{\dimPred,\rm{max}}^{\ubm}$ (resp. $\e_{\dimPred,\rm{max}}^{\nug\nug}$) in the left (resp. right) panel in order to compare the displacement (resp. normal stress) error between the mixed and Nitsche approaches.}
Considering the mixed approach, we notice a stagnation for all errors after a certain number of primal modes $N$. For the primal error, there is a clear decrease of the error as a function of the dimension of the dual basis $R$. However, for the dual error, although the error decreases, we observe some oscillations for certain values of the dimension of the primal basis $N$ (compare the errors for $R=30$ and $R=40$). This can be explained by the fact that the dual basis obtained with the mCPG algorithm is not orthonormal (owing to the sign constraints on the Lagrange multiplier, since an orthonormalization process cannot be performed). Moreover, the higher the number of vectors in the dual basis, the more the orthogonality property is lost, and this fact introduces noise in the reduced model. In terms of accuracy, we observe that we have a better approximation for the primal variable (error of the order of $10^{-4}\SI{}{\metre}$) than for the dual variable (error of the order  of $10^{-2}\SI{}{\pascal}$) in the mixed formulation. These results illustrate the fact that it is very difficult to reduce the dual space, and therefore further motivate the use of purely primal methods like Nitsche's method in the RBM framework applied to contact problems. {Finally, comparing the mixed and Nitsche formulations, we observe that the displacement error for the latter is better than for the former for the three values of the dimension of the dual reduced cone. Instead, comparing the error on the normal stress, we observe that the error for the former is slightly better when $R=40$, but the dimension of the primal reduced space is much larger owing to the need to ensure inf-sup stability.}

\subsection{Friction case}\label{subsection:4:FC}
Let us now consider the Tresca frictional Hertz contact problem. We choose a threshold $s=0.1$. The parameter $\gamma$ and the mesh size $h$ are the same as for the frictionless case (see Section~\ref{subsubsection:4:N}).

\begin{figure}[!h]
	\centering
	\begin{subfigure}{0.32\textwidth}
		\includegraphics[width=0.8\textwidth]{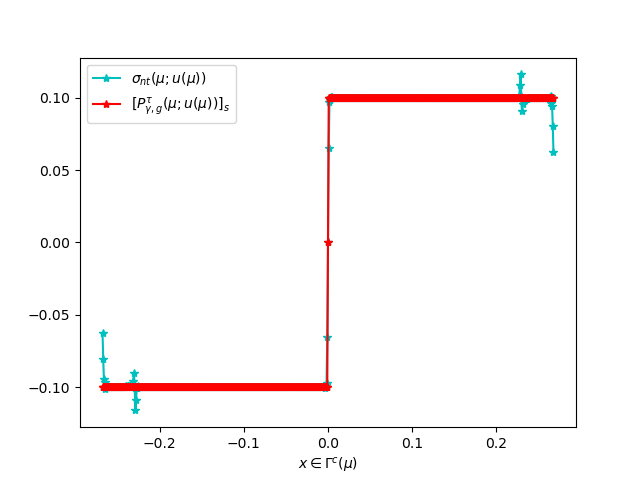}
		%\caption{$\sigman(\ubm(\mu))$.}
		\label{fig:4:HertzLNT:sigma07}
	\end{subfigure}
	\hfill
	\begin{subfigure}{0.32\textwidth}
		\includegraphics[width=0.8\textwidth]{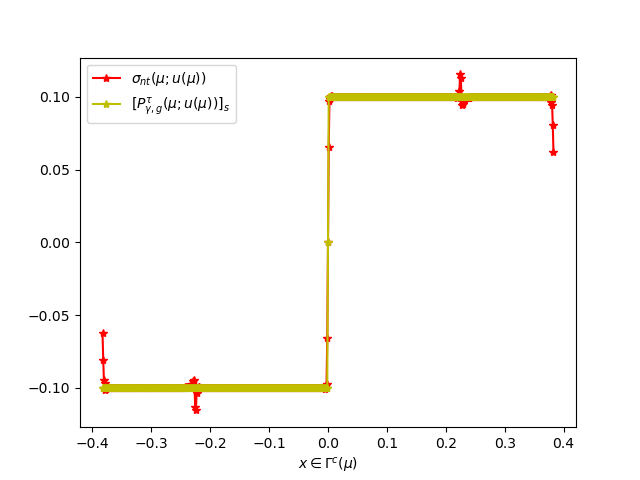}
		%\caption{$\un(\mu) - g(\mu)$.}
		\label{fig:4:HertzLNT:sigma10}
	\end{subfigure}
	\hfill
	\begin{subfigure}{0.32\textwidth}
		\includegraphics[width=0.8\textwidth]{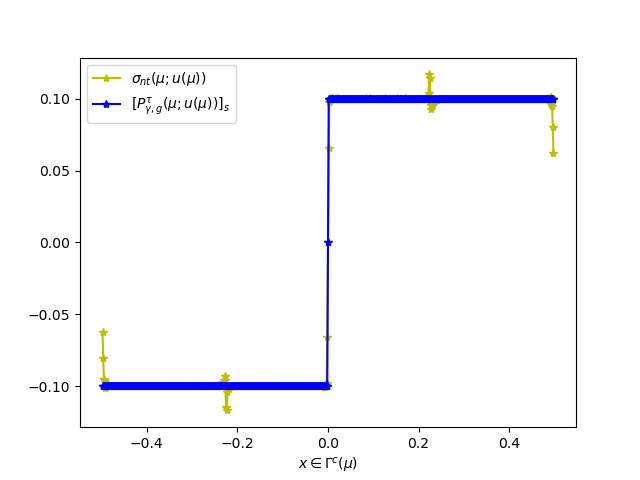}
		%\caption{$\un(\mu) - g(\mu)$.}
		\label{fig:4:HertzLNT:sigma13}
	\end{subfigure}
	
	\caption{{Frictional Hertz test case; Nitsche's method (HF). $\sigmat(\ubm(\mu))$ and $\projconvb{\Pgamt(\mu;\ubm(\mu))}{s}$. Column~$1$: $\mu=0.7\SI{}{\metre}$. Column~$2$: $\mu=1.0\SI{}{\metre}$. Column~$3$: $\mu=1.3\SI{}{\metre}$.}}
	\label{fig:4:HertzLNT:CU-AL}
\end{figure}

Figure~\ref{fig:4:HertzLNT:CU-AL} displays the superposition of the tangential stress $\sigmat(\ubm(\mu))$ and its Alart--Curnier reformulation  $\projconvb{\Pgamt(\mu;\ubm(\mu))}{s}$ as a function of the abscissa along $\gammaac(\mu)$ for $\mu = 0.7\SI{}{\metre}$ (Column~$1$), $\mu = 1\SI{}{\metre}$ (Column~$2$), and $\mu = 1.3\SI{}{\metre}$ (Column~$3$). We observe that the tangential stress matches very well with its counterpart $\projconvb{\Pgamt(\mu;\ubm(\mu))}{s}$ resulting from the Alart--Curnier reformulation. However, we observe some oscillations on $\sigmat(\ubm(\mu))$ at the end of the effective contact zone (contact/non-contact transition zone) and at the end of the potential contact zone. We also notice that the tangential stress is not zero outside the effective contact zone, which is consistent with a known (undesirable) feature of Tresca's model, namely to predict friction without contact. 
To have a better look at the  friction conditions, we display in Table~\ref{table:4:HertzLNTC:signo} the relative errors on the Alart--Curnier reformulation of the Tresca friction conditions defined as follows:
\begin{equation}
e_{\rm{AC}}^{\nug\taug,\rm{T}}(\mu):= \frac{\norm{\sigmat(\ubm(\mu)) - \projconvb{\Pgamt(\mu;\ubm(\mu))}{s}}{\ell^{2}(\gammaac(\mu))}}{\norm{\sigman(\ubm(\mu))}{\ell^{2}(\gammaac(\mu))}}.
\end{equation}
For the three values of the parameter $\mu \in \{0.7,1,1.3\}(\SI{}{\metre})$, we consider three values of the mesh size $h$, namely $h=5\SI{}{\milli\metre}$ (coarse), $h=2.5\SI{}{\milli\metre}$ (medium) and $h=1.25\SI{}{\milli\metre}$ (fine). We notice that the relative error $e_{\rm{AC}}^{\nug\taug,\rm{T}}(\mu)$ is smaller than $6\%$ for the three parameters values and for the three mesh sizes. Moreover, we see that the errors decrease when $h$ decreases, thereby indicating the convergence of the approximation. More precisely, we observe a convergence of order $\frac{1}{2}$ for the approximation of $\sigmat(\ubm(\mu))$. Notice that this order is different from the one observed in the frictionless case (order $1$). Thus, we can say that the Tresca friction conditions are globally satisfied with a good accuracy as for the Signorini contact conditions although they are not strictly enforced. 
\begin{table}[!h]
	\centering
	\renewcommand{\arraystretch}{1.5}
	\begin{tabular}{|c|*{9}{>{\centering\arraybackslash}m{1.1cm}|} }
		\hline
		$\mu(\SI{}{\metre})$ & \multicolumn{3}{|c|}{$0.7$}&  \multicolumn{3}{|c|}{$1$}&  \multicolumn{3}{|c|}{$1.3$} \\
		\hline
		$h(\SI{}{\milli\metre})$&$5$ & $2.5$& $1.25$ &$5$ & $2.5$& $1.25$ &$5$ & $2.5$& $1.25$ \\
		\hline
		$e_{\rm{AC}}^{\nug\taug,\rm{T}}(\%)$ &$ 5.49$& $ 3.33$& $ 2.36$& $ 5.5$ &$ 3.35$& $2.33 $& $ 5.76$& $3.4 $&$ 2.37$ \\
		\hline
	\end{tabular}
	\caption{{Frictional Hertz test case, Nitsche's method (HF): Relative errors $e_{\rm{AC}}^{\nug\taug,\rm{T}}(\mu)$ for $\mu \in \{0.7,1.0,1.3\}(\SI{}{\metre})$ and the mesh sizes $h \in \{5,2.5,1.25\}(\SI{}{\milli\metre})$.}}
	\label{table:4:HertzLNTC:signo}
\end{table}

\begin{figure}[!h]
	\centering
	\begin{minipage}[l]{.45\linewidth}
		\includegraphics[width=7cm]{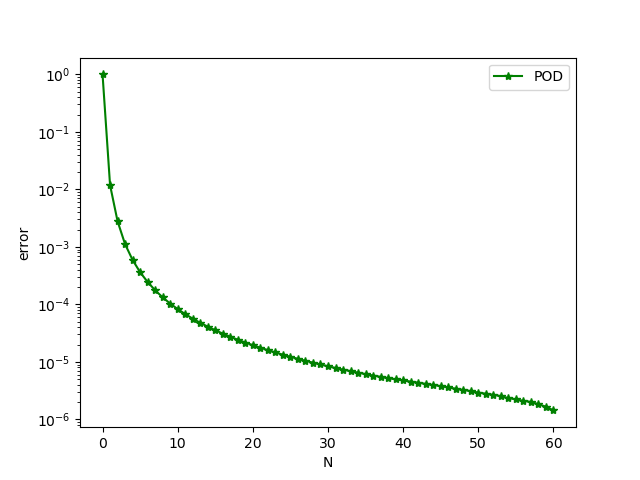}
	\end{minipage} 
	\caption {{Frictional Hertz test case, Nitsche's method: Relative POD projection error $e_{\texttt{POD}}(N)$ as function of $\dimPred$, for $|\discrd{D}{train}|=61$.}}
\label{fig:4:HertzLNT:POD_error}
\end{figure}

{Figure~\ref{fig:4:HertzLNT:POD_error} shows the relative projection error $e_{\texttt{POD}}(N)$ produced by the POD algorithm (see~\eqref{equa:4:PODerror}) as a function of the number of vectors composing the reduced basis. We notice that the projection error decreases sufficiently fast so that indeed the linear spaces generated by the snapshots can be approximated by small-dimensional subspaces. We also observe a fast decrease of the POD error for the first $15$ modes before a slower decrease occurs at relative error levels between $10^{-5}$ and $10^{-6}$.}

\begin{figure}[!h]
	\centering
	\begin{minipage}[l]{.45\linewidth}
		\includegraphics[width=7cm]{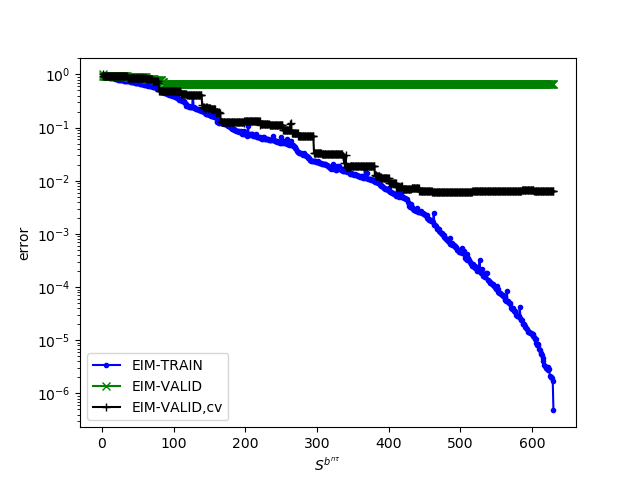}
	\end{minipage} 
	\begin{minipage}[l]{.45\linewidth}
		\includegraphics[width=7cm]{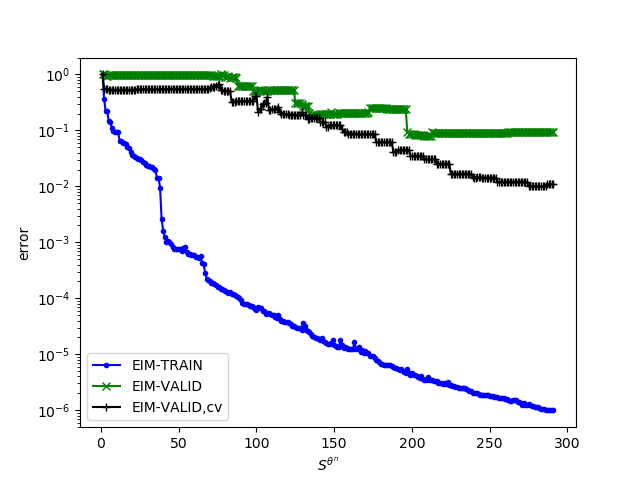}
	\end{minipage}
		\caption {{Frictional Hertz test case, Nitsche's method: Relative EIM approximation errors as a function of the rank $\eimcar{b^{\nug\taug}}$ or $\eimcar{\thetan}$ of the approximation. Left: $e_{\texttt{EIM}}^{b^{\nug\taug}}(\eimcar{b^{\nug\taug}},\discrd{D}{train})$, $e_{\texttt{EIM}}^{b^{\nug\taug}}(\eimcar{b^{\nug\taug}},\discrd{D}{valid})$ and $e_{\texttt{EIM}}^{b^{\nug\taug},\rm{cv}}(\eimcar{b^{\nug\taug}})$; Right: $e_{\texttt{EIM}}^{\thetan}(\eimcar{\thetan},\discrd{D}{train})$, $e_{\texttt{EIM}}^{\thetan}(\eimcar{\thetan},\discrd{D}{valid})$ and $e_{\texttt{EIM}}^{\thetan,\rm{cv}}(\eimcar{\thetan})$.}}
	\label{fig:4:HertzLNT:eim_error}
\end{figure}

Let us discuss the EIM approximation for which the training sets $\eimtr{b^{\nug\taug}}$, $\eimtr{\thetan}$ and $\eimtr{\thetat}$ are of cardinality $898$. We fix a tolerance $\tole{EIM}:=10^{-6}$. With this choice, we obtain $\eimcar{b^{\nug\taug}}=630\ll \dimP\times \dimP$,  $\eimcar{\thetan}=291 \ll \dimP$ and  $\eimcar{\thetat}=3 \ll \dimP$. For the same reason as for the frictionless case, we do not need to perform an EIM decomposition on $\Agamm(\mu)$ and on $\F(\mu)$.  Figure~\ref{fig:4:HertzLNT:eim_error} shows the relative EIM interpolation errors for the tangent matrix $\Bgamm(\mu,\ubmka{k}(\mu))$ (left panel) and the residual vector $\Rhonm(\mu,\ubmka{k}(\mu))$ (right panel) as a function of the rank $\eimcar{b^{\nug\taug}}$ or $\eimcar{\thetan}$, \ie, we plot $e_{\texttt{EIM}}^{\thetan}(\eimcar{\thetan},\discrd{D}{*})$ defined in~\eqref{equa:4:erroreimthetan} and 
\begin{equation}
	{e_{\texttt{EIM}}^{b^{\nug\taug}}(\eimcar{b^{\nug\taug}},\discrd{D}{*}):= \frac{\max\limits_{\mu \in \discrd{D}{*}}\  \max\limits_{k \in \indset{1}{\kc(\mu)}} \norm{\Bgamm(\mu,\ubmka{k}(\mu)) - \eimop{b^{\nug\taug}}(\mu,k)}{\ell^{\infty}(ij)}}{\max\limits_{\mu \in \discrd{D}{*}}\ \max\limits_{k \in \indset{1}{\kc(\mu)}} \norm{\Bgamm(\mu,\ubmka{k}(\mu))}{\ell^{\infty}(ij)}},}
\end{equation}
{with $\mathcal{D}_{*}$ either equal to $\discrd{D}{train}$ or to $\discrd{D}{valid}$.
For the validation set $\discrd{D}{valid}$, we additionally plot the relative EIM approximation errors  
\begin{equation}
e_{\texttt{EIM}}^{b^{\nug\taug},\rm{cv}}(\eimcar{b^{\nug\taug}}):= \frac{\max\limits_{\mu \in \discrd{D}{valid}}  \norm{\Bgamm(\mu,\ubmka{\rm{cv}}(\mu)) - \eimop{b^{\nug\taug}}(\mu,\kc(\mu))}{\ell^{\infty}(ij)}}{\max\limits_{\mu \in \discrd{D}{valid}} \norm{\Bgamm(\mu,\ubmka{\rm{cv}}(\mu))}{\ell^{\infty}(ij)}},
\end{equation}
and $e_{\texttt{EIM}}^{\thetan,\rm{cv}}(\eimcar{\thetan})$ defined in~\eqref{equa:4:erroreimthetanconv}.
Notice that for Tresca friction, we perform the EIM on the tangent matrix $\Bgamm(\mu;\ubmka{k}(\mu))$ instead of performing it separately on the tangent matrices $\Bgamnm(\mu;\ubmka{k}(\mu))$ and $\Bgamtm(\mu;\ubmka{k}(\mu))$. We see that the result is close to the one obtained for the frictionless case because the contribution of $\Bgamtm(\mu;\ubmka{k}(\mu))$ is negligible. Indeed, as can be seen in Figure~\ref{fig:4:HertzLNT:CU-AL} (Row~$1$), the Alart--Curnier reformulation of the tangential stress is equal to $\pm s$ almost everywhere (except at $3$ nodes) which results in the nullity of $\Gdif{s}(\Pgamt(\mu; \ubmka{k}(\mu)))$ almost everywhere. For the same reason, the dependence of the residual vector $\Rhotm(\mu,\ubmka{k}(\mu))$ on $(\mu,k)$ is almost of rank one and therefore its EIM approximation is very simple ($\eimcar{\thetat}=3$); hence, its relative EIM approximation error is not reported. Altogether, the behaviour of the residual vector $\Rhonm(\mu,\ubmka{k}(\mu))$ remains similar to that observed in the frictionless case.} 

\begin{figure}[!h]
	\centering
	\begin{minipage}[l]{.45\linewidth}
		\includegraphics[width=7cm]{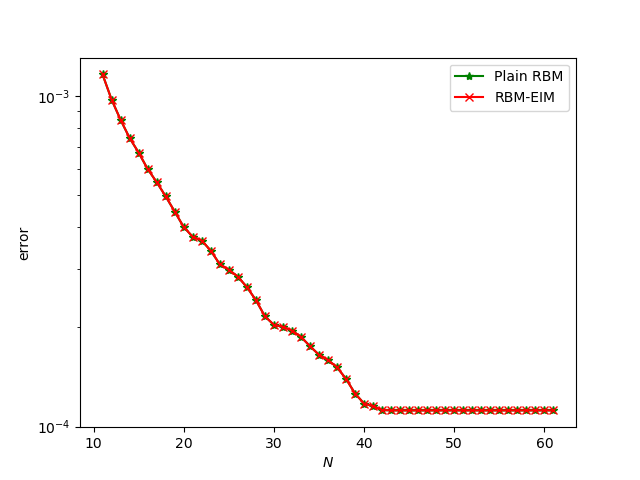}
	\end{minipage} 
	\begin{minipage}[l]{.45\linewidth}
		\includegraphics[width=7.5cm]{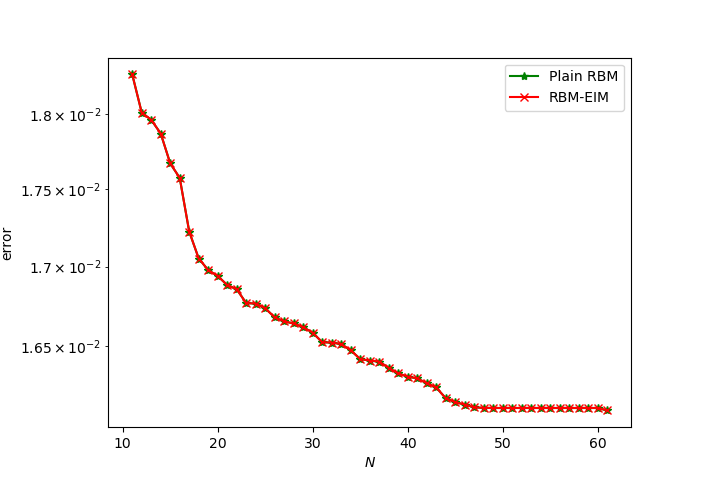}
	\end{minipage}
	\caption {{Frictional Hertz test case, Nitsche's method: RBM approximation errors  for $|\discrd{D}{train}|=61$. Left: $\e_{\dimPred,\rm{max}}^{\ubm}$; Right: $\e_{\dimPred,\rm{max}}^{\nug\taug}$.}}
	\label{fig:4:HertzLNT:RB_error}
\end{figure}

We denote by $\erb{\nug\taug}(\mu)$ the relative RB approximation error on the  tangential stress defined as
\begin{equation}
\erb{\nug\taug}(\mu):=\frac{\norm{\sigmat(\ubmka{\rm{cv}}(\mu)) - \sigmat(\ubmredka{\rm{cv}}(\mu))}{\ell^2(\gammaac(\mu))}}{\norm{\sigmat(\ubmka{\rm{cv}}(\mu)) }{\ell^2(\gammaac(\mu))}},
\end{equation}
and introduce the relative error measure $\e_{\dimPred,\rm{max}}^{\nug\taug}$  defined as
\begin{equation}
	\e_{\dimPred,\rm{max}}^{\nug\taug}:= \max\limits_{ \mu\in \discrd{D}{valid}} \erb{\nug\taug}(\mu).
\end{equation}
Figure~\ref{fig:4:HertzLNT:RB_error} displays for Tresca's friction the relative errors $\e_{\dimPred,\rm{max}}^{\ubm}$ (left panel) and $\e_{\dimPred,\rm{max}}^{\nug\taug}$ (right panel) as a function of the number of vectors composing the reduced basis. As for the frictionless case, we only consider RB dimensions $N$ larger than $10$ and the validation set  $\discrd{D}{valid}$ is restricted to the interval $[0.7,1.18](\SI{}{\metre})$. In Figure~\ref{fig:4:HertzLNT:RB_error}, we superpose the plain RBM approximation error (without any EIM) and the RBM-EIM approximation error. As for the frictionless case, we observe that similar errors are obtained for plain RBM and RBM-EIM. This again confirms the satisfactory quality of the EIM approximations at convergence. 

\subsection{Toward Coulomb friction} \label{sec:Coulomb}

A more realistic model for friction is given by Coulomb conditions which read as follows: 
\begin{equation}
\left\{
\begin{alignedat}{2}
\|\sigmat(\ubm(\mu))\| \leq \nuf |\sigman(\ubm(\mu))|,&\quad \mbox{if }  \ut(\mu) =\zerobm,\\
\sigmat(\ubm(\mu)) = -\nuf |\sigman(\ubm(\mu))| \frac{\ut(\mu)}{\|\ut(\mu)\|}, &\quad \mbox{otherwise},
\end{alignedat}
\right.
\label{equa:4:coul}
\end{equation}
where  $\nuf >0$ is a given nondimensional coefficient (which can be taken to be constant for simplicity).
A rather well-established approach to solve the Coulomb frictional problem is to use a fixed-point method on the Tresca frictional problem. Specifically, one introduces the mapping $\phibm(\mu;\cdot): \fespacebm{V}(\mu) \rightarrow \fespacebm{V}(\mu)$ defined for all $\vbm \in \fespacebm{V}(\mu)$ by requiring that $\phibm(\mu;\vbm)$ solves~\eqref{equa:4:NWF} with the threshold $s:=\snc(\mu;\vbm)$ with $\snc(\mu;\vbm):= \nuf\big|\projconv{\Pgamn(\mu;\vbm)}{-}\big|$. It is shown in~\cite[Thm~$4.3$]{Chouly-2022} that, under certain conditions, the mapping $\phibm(\mu;\cdot)$ is contractive. Thus, the Coulomb frictional problem can be solved by means of a nested loop, where the outer iteration index, say $n$, refers to the fixed-point iteration on $\phibm(\mu;\cdot)$, and the inner iteration index, say $k$, refers as in~\eqref{equa:4:NWF} to the iterative solution of the Tresca frictional problem. 

\begin{figure}[!h]
	\centering	
	\begin{subfigure}{0.32\textwidth}
		\includegraphics[width=0.8\textwidth]{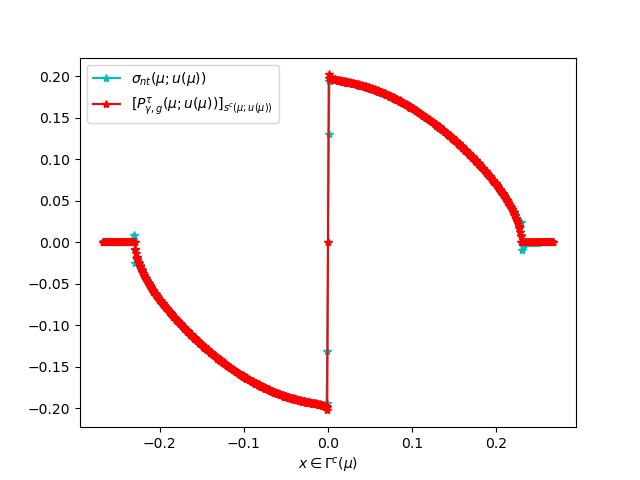}
		%\caption{$\sigman(\ubm(\mu))$.}
		\label{fig:4:HertzLNC:sigma07}
	\end{subfigure}
	\hfill
	\begin{subfigure}{0.32\textwidth}
		\includegraphics[width=0.8\textwidth]{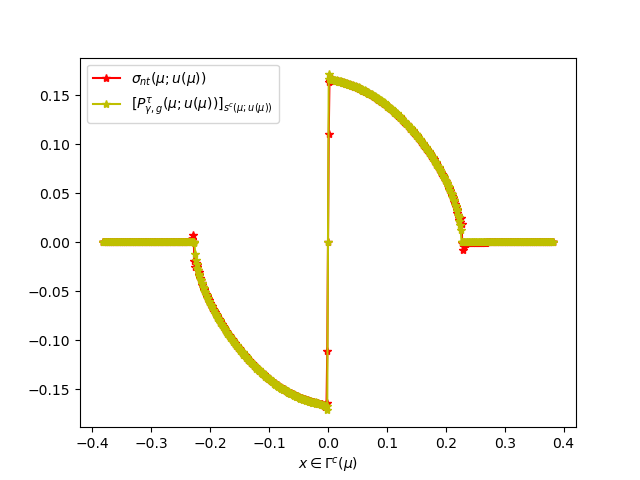}
		%\caption{$\un(\mu) - g(\mu)$.}
		\label{fig:4:HertzLNC:sigma10}
	\end{subfigure}
	\hfill
	\begin{subfigure}{0.32\textwidth}
		\includegraphics[width=0.8\textwidth]{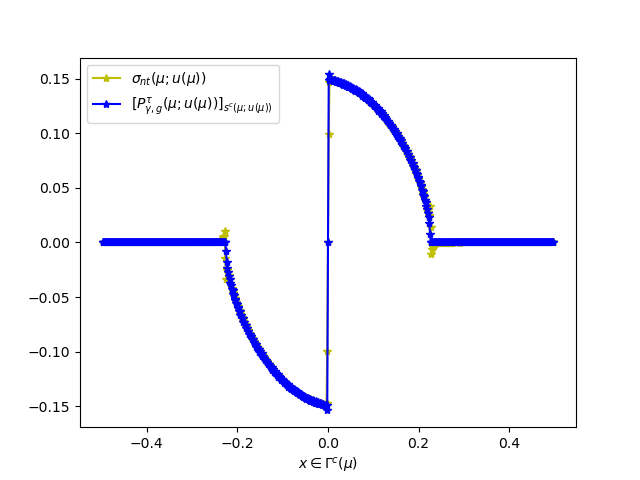}
		%\caption{$\un(\mu) - g(\mu)$.}
		\label{fig:4:HertzLNC:sigma13}
	\end{subfigure}
	
	\caption{{Frictional Hertz test case with Coulomb friction law, Nitsche's method (HF).  $\sigmat(\ubm(\mu))$ and $\projconvb{\Pgamt(\mu;\ubm(\mu))}{\snc(\mu;\ubm(\mu))}$. Column~$1$: $\mu=0.7\SI{}{\metre}$. Column~$2$: $\mu=1.0\SI{}{\metre}$. Column~$3$: $\mu=1.3\SI{}{\metre}$.}}
	\label{fig:4:HertzLNC:CU-AL}
\end{figure}

HF solutions for the Hertz test case obtained using the above nested loop can be obtained to generate snapshots for the parameter values from the training set. Three snapshots are illustrated in Figure~\ref{fig:4:HertzLNC:CU-AL} (compare with Figure~\ref{fig:4:HertzLNT:CU-AL}). 
Moreover, the set of snapshots can be compressed by using POD, leading to similar results to those obtained for Tresca friction (see Figure~\ref{fig:4:HertzLNT:POD_error}). The main challenge within the current approach lies in the realization of the EIM since the nested iterative loop now requires to separate the dependencies on the triple $(\mu,n,k)$. To overcome this difficulty, one possibility is to consider only converged solutions for the inner iteration (index $k$) when computing the EIM decompositions of the tangent matrix and the residual. However, these decompositions turn out to be, so far, rather inaccurate at the early stages of the iterative procedure, thereby hampering convergence. This difficulty will be further investigated in future work.

\bibliographystyle{plain} 
\bibliography{biblio.bib} 

\begin{thebibliography}{10}

\bibitem{Baillet-2006}
L.~Baillet and T.~Sassi.
\newblock Mixed finite element methods for the {S}ignorini problem with
  friction.
\newblock {\em Numer. Methods Partial Differential Equations},
  22(6):1489--1508, 2006.

\bibitem{Balajewicz-2016}
M.~Balajewicz, D.~Amsallem, and C.~Farhat.
\newblock Projection-based model reduction for contact problems.
\newblock {\em Internat. J. Numer. Methods Engrg.}, 106(8):644--663, 2016.

\bibitem{Barrault-2004}
M.~Barrault, Y.~Maday, N.~C. Nguyen, and A.~T. Patera.
\newblock An `empirical interpolation' method: application to efficient
  reduced-basis discretization of partial differential equations.
\newblock {\em C. R. Math. Acad. Sci. Paris}, 339(9):667--672, 2004.

\bibitem{Benaceur-2020}
A.~Benaceur, A.~Ern, and V.~Ehrlacher.
\newblock A reduced basis method for parametrized variational inequalities
  applied to contact mechanics.
\newblock {\em Internat. J. Numer. Methods Engrg.}, 121(6):1170--1197, 2020.

\bibitem{Buffa-2012}
A.~Buffa, Y.~Maday, A.~T. Patera, C.~Prud'homme, and G.~Turinici.
\newblock {\it {A} priori} convergence of the greedy algorithm for the
  parametrized reduced basis method.
\newblock {\em ESAIM Math. Model. Numer. Anal.}, 46(3):595--603, 2012.

\bibitem{Chouly-2014}
F.~Chouly.
\newblock An adaptation of {N}itsche's method to the {T}resca friction problem.
\newblock {\em J. Math. Anal. Appl.}, 411(1):329--339, 2014.

\bibitem{Chouly-2020}
F.~Chouly, A.~Ern, and N.~Pignet.
\newblock A hybrid high-order discretization combined with {N}itsche's method
  for contact and {T}resca friction in small strain elasticity.
\newblock {\em SIAM J. Sci. Comput.}, 42(4):A2300--A2324, 2020.

\bibitem{Chouly-2017}
F.~Chouly, M.~Fabre, P.~Hild, R.~Mlika, J.~Pousin, and Y.~Renard.
\newblock An overview of recent results on {N}itsche's method for contact
  problems.
\newblock In {\em Geometrically unfitted finite element methods and
  applications}, volume 121 of {\em Lect. Notes Comput. Sci. Eng.}, pages
  93--141. Springer, Cham, 2017.

\bibitem{Chouly-2013}
F.~Chouly and P.~Hild.
\newblock A {N}itsche-based method for unilateral contact problems: numerical
  analysis.
\newblock {\em SIAM J. Numer. Anal.}, 51(2):1295--1307, 2013.

\bibitem{Chouly-2022}
F.~Chouly, P.~Hild, V.~Lleras, and Y.~Renard.
\newblock Nitsche method for contact with {C}oulomb friction: existence results
  for the static and dynamic finite element formulations.
\newblock {\em J. Comput. Appl. Math.}, 416, 2022.

\bibitem{Chouly-2015}
F.~Chouly, P.~Hild, and Y.~Renard.
\newblock Symmetric and non-symmetric variants of {N}itsche's method for
  contact problems in elasticity: theory and numerical experiments.
\newblock {\em Math. Comp.}, 84(293):1089--1112, 2015.

\bibitem{Curnier-1988}
A.~Curnier and P.~Alart.
\newblock A generalized {N}ewton method for contact problems with friction.
\newblock {\em J. M\'{e}c. Th\'{e}or. Appl.}, 7(suppl. 1):67--82, 1988.

\bibitem{Duvaut-1972}
G.~Duvaut and J.-L. Lions.
\newblock {\em Les in\'{e}quations en m\'{e}canique et en physique}.
\newblock Travaux et Recherches Math\'{e}matiques, No. 21. Dunod, Paris, 1972.

\bibitem{Fauque-2018}
J.~Fauque, I.~Rami\`ere, and D.~Ryckelynck.
\newblock Hybrid hyper-reduced modeling for contact mechanics problems.
\newblock {\em Internat. J. Numer. Methods Engrg.}, 115(1):117--139, 2018.

\bibitem{Fichera-1964}
G.~Fichera.
\newblock Problemi elastostatici con vincoli unilaterali: {I}l problema di
  {S}ignorini con ambigue condizioni al contorno.
\newblock {\em Atti Accad. Naz. Lincei Mem. Cl. Sci. Fis. Mat. Natur. Sez. Ia
  (8)}, 7:91--140, 1963/64.

\bibitem{Fortin-1983}
M.~Fortin and R.~Glowinski.
\newblock {\em Augmented {L}agrangian methods}, volume~15 of {\em Studies in
  Mathematics and its Applications}.
\newblock North-Holland Publishing Co., Amsterdam, 1983.
\newblock Applications to the numerical solution of boundary value problems,
  Translated from the French by B. Hunt and D. C. Spicer.

\bibitem{Haasdonk-2013}
B.~Haasdonk.
\newblock Convergence rates of the {POD}-greedy method.
\newblock {\em ESAIM Math. Model. Numer. Anal.}, 47(3):859--873, 2013.

\bibitem{Haasdonk-2012}
B.~Haasdonk, J.~Salomon, and B.~Wohlmuth.
\newblock A reduced basis method for parametrized variational inequalities.
\newblock {\em SIAM J. Numer. Anal.}, 50(5):2656--2676, 2012.

\bibitem{hesthaven2016certified}
J.~S. Hesthaven, G.~Rozza, and B.~Stamm.
\newblock {\em Certified reduced basis methods for parametrized partial
  differential equations}, volume 590.
\newblock Springer, 2016.

\bibitem{Johnson-1987}
K.~L. Johnson.
\newblock {\em Contact mechanics}.
\newblock Cambridge University Press, 1987.

\bibitem{Kikuchi-1988}
N.~Kikuchi and J.~T. Oden.
\newblock {\em Contact {P}roblems in {E}lasticity: {A} {S}tudy of {V}ariational
  {I}nequalities and {F}inite {E}lement {M}ethods}, volume~8 of {\em SIAM
  Studies in Applied Mathematics}.
\newblock Society for Industrial and Applied Mathematics (SIAM), Philadelphia,
  PA, 1988.

\bibitem{kollepara2023limitations}
K.~S. Kollepara, J.~M. Navarro-Jim{\'e}nez, Y.~Le~Guennec, L.~Silva, and J.~V.
  Aguado.
\newblock On the limitations of low-rank approximations in contact mechanics
  problems.
\newblock {\em Internat. J. Numer. Methods Engrg.}, 124(1):217--234, 2023.

\bibitem{Kunisch-2001}
K.~Kunisch and S.~Volkwein.
\newblock Galerkin proper orthogonal decomposition methods for parabolic
  problems.
\newblock {\em Numer. Math.}, 90(1):117--148, 2001.

\bibitem{LeBRFD:23}
S.~Le~Berre, I.~Rami{\`e}re, J.~Fauque, and D.~Ryckelynck.
\newblock Condition number and clustering-based efficiency improvement of
  reduced-order solvers for contact problems using {L}agrange multipliers.
\newblock {\em Mathematics}, 10:1495--1520, 2022.

\bibitem{Maday-2009}
Y.~Maday, N.~C. Nguyen, A.~T. Patera, and G.~S.~H. Pau.
\newblock A general multipurpose interpolation procedure: the magic points.
\newblock {\em Commun. Pure Appl. Anal.}, 8(1):383--404, 2009.

\bibitem{Mlika-2017}
R.~Mlika, Y.~Renard, and F.~Chouly.
\newblock An unbiased {N}itsche's formulation of large deformation frictional
  contact and self-contact.
\newblock {\em Comput. Methods Appl. Mech. Engrg.}, 325:265--288, 2017.

\bibitem{Niakh-2022}
I.~Niakh, G.~Drouet, V.~Ehrlacher, and A.~Ern.
\newblock {Stable model reduction for linear variational inequalities with
  parameter-dependent constraints}.
\newblock M2AN, to appear,
  \texttt{https://hal.archives-ouvertes.fr/hal-03611982}, March 2022.

\bibitem{Nitsche-1971}
J.~Nitsche.
\newblock \"{U}ber ein {V}ariationsprinzip zur {L}\"{o}sung von
  {D}irichlet-{P}roblemen bei {V}erwendung von {T}eilr\"{a}umen, die keinen
  {R}andbedingungen unterworfen sind.
\newblock {\em Abh. Math. Sem. Univ. Hamburg}, 36:9--15, 1971.

\bibitem{prud2002reliable}
C.~Prud’Homme, D.~V. Rovas, K.~Veroy, L.~Machiels, Y.~Maday, A.~T. Patera,
  and G.~Turinici.
\newblock Reliable real-time solution of parametrized partial differential
  equations: Reduced-basis output bound methods.
\newblock {\em J. Fluids Eng.}, 124(1):70--80, 2002.

\bibitem{quarteroni2015reduced}
A.~Quarteroni, A.~Manzoni, and F.~Negri.
\newblock {\em Reduced basis methods for partial differential equations: an
  introduction}, volume~92.
\newblock Springer, 2015.

\bibitem{Renard-2020}
Y.~Renard and K.~Poulios.
\newblock {GetFEM: Automated FE modeling of multiphysics problems based on a
  generic weak form language}.
\newblock \texttt{https://hal.archives-ouvertes.fr/hal-02532422}, 2020.

\bibitem{Sofonea-2012}
M.~Sofonea and A.~Matei.
\newblock {\em Mathematical models in contact mechanics}, volume 398.
\newblock Cambridge University Press, 2012.

\bibitem{Stampacchia-1964}
G.~Stampacchia.
\newblock Formes bilin\'{e}aires coercitives sur les ensembles convexes.
\newblock {\em C. R. Acad. Sci. Paris}, 258:4413--4416, 1964.

\bibitem{Wriggers-2006}
P.~Wriggers.
\newblock {\em Computational Contact Mechanics}, volume~2.
\newblock Springer Berlin, Heidelberg, 2006.

\bibitem{zeka2022preliminary}
D.~Zeka, P.-A. Guidault, D.~N{\'e}ron, M.~Guiton, and G.~Ench{\'e}ry.
\newblock Preliminary study for the simulation of wire ropes using a model
  reduction approach suitable for multiple contacts.
\newblock In {\em 25{\`e}me Congr{\`e}s Fran{\c{c}}ais de M{\'e}canique}, 2022.

\end{thebibliography}
\end{document}